\documentclass[12pt]{amsart}
\usepackage{amscd}

\newtheorem{thm}{Theorem}[section]
\newtheorem{con}[thm]{Conjecture}

\theoremstyle{definition}

\def \ph{\varphi}

\def\aut{\operatorname{Aut}}
\def\laut{\operatorname{LAut}}
\def\gaut{\operatorname{GAut}}
\def\autf{\operatorname{Aut_f}}

\def \diag{\operatorname {diag}}

\def \refeq#1{equation (\ref{#1})}
\def \ra{\rightarrow}

\def \hom{\mbox{\rm Hom}}

\def \mplus{+\cdots+}

\def \Z{\mbox{$\mathbb Z$}}

\def\br#1#2{\lbrack#1,#2\rbrack}

\def\zt{\mbox{$\Z_2$}}

\def\ad{\operatorname{ad}}

\def\inv{^{-1}}
\def\d{d}

\def\D{D}

\def\coder{\operatorname{Coder}}

\def\linf{\mbox{$L_\infty$}}
\def\and{\mbox{ \rm and }}

\def\s#1{(-1)^{#1}}

\def\phd#1#2{\ph^{#1}_{#2}}

\def\thmref#1{Theorem (\ref{#1})}
\def\inv{^{-1}}
\def\psd#1#2{\psi^{#1}_{#2}}
\def\dstar{\mbox{$d_*$}}

\def\dsharp{\mbox{$d_{\sharp}$}}
\def\dzeroa{\mbox{$d_{0,a}$}}

\def\dt#1{\mbox{$d_{#1}$}}
\def\ho{\text{ho}}
\def\su{\bigl(\tfrac su\bigr)}
\def\pfrac#1#2{\bigl(\tfrac{#1}{#2}\bigr)}
\def\dfin{d_{\text{fin}}}

\author{Alice Fialowski}
\address{E\"otv\"os Lor\'and University\\ Budapest, Hungary}
\email{fialowsk@cs.elte.hu}
\author{Michael Penkava}
\address{University of Wisconsin\\ Eau Claire, WI 54702-4004}
\email{penkavmr@uwec.edu}
\thanks{The research
of the authors was supported by grants MTA-OTKA-NSF 38453, OTKA T043641
and T043641 and by grants from the University of Wisconsin-Eau Claire}
\title[Extensions of $2|1$ Dimensional Algebras]{Extensions of \linf\
algebras of two even and one odd dimension}

\begin{document}
\setlength{\multlinegap}{0pt}

\begin{abstract} In this article we study extensions of  $\Z_2$-graded
\linf\ algebras on a vector space of two even and one odd dimension.
  In particular, we
determine all extensions of a super Lie algebra as an \linf\ algebra.
Our convention on the parities is the opposite of the usual one,
because we define our structures on the symmetric coalgebra of the
parity reversion of a space. \end{abstract}
\date\today \maketitle


\section{Introduction}
The problem of classifying formal deformations of algebras has a long and interesting history.
The formality conjecture, proved by Maxim Kontsevich, is an example of a problem where the
existence of an extension of an infinitesimal deformation of an algebra to a formal deformation
has important applications in both mathematics and mathematical physics.  In addition, the classification
of all such extensions is an important question that is not so easy to solve.  A miniversal deformation
of an algebra contains the information necessary to construct all possible formal extensions of an
infinitesimal deformation. In a series of articles, the authors have been studying
examples of miniversal
deformations of low dimensional \linf\ algebras, with the aim of developing a constructive approach to
the process.

 In \cite{bfp1},
codifferentials of degree one and two on a $2|1$-dimensional vector space were studied and miniversal
deformations of  degree two codifferentials were constructed. Since degree
two codifferentials are \zt-graded Lie algebra structures, the miniversal
deformations describe the deformations of these Lie algebra structures into
more general \linf\ algebras.
In the current article, which is a continuation of the study of \linf\ structures
on a $2|1$ dimensional space, we focus on another manner in which \linf\
algebras are constructed from \zt-graded Lie algebras, by extending the
codifferential to an \linf\ algebra by adding higher order terms.

The two ideas, extensions of \linf\ algebras, and deformations of \linf\ algebras
are closely related. In fact, an extension of a degree $n$ codifferential to a
\linf\ algebra can be thought of as a special case of a deformation, where the deformation parameters
are assigned fixed values. However, our point of view is different in this paper.
We will be interested in classifying all extensions of a fixed codifferential
of degree 2 to an \linf\ algebra up to equivalence, where equivalent extensions are
determined by automorphisms of the symmetric coalgebra of the underlying \zt-graded
vector space.

This is a far more difficult problem than simply determining all the
\linf\ algebra structures. In fact, for a $2|1$-dimensional vector space, it is trivial
to give a list of all \linf\ structures; there are in fact two distinct \emph{kinds}
of odd cochains, and an odd coderivation is a codifferential precisely when it contains
only terms of one or the other kind of odd cochains.

In the physics literature, \linf\ algebras are usually referred to as
strongly homotopy Lie algebras, or sh-Lie algebras. These algebras first appeared in
\cite{ss}, and have been studied in mathematics (\cite{sta3}, \cite{hsch}) as well as in
mathematical physics (\cite{ls}, \cite{z}, \cite{aksz},
\cite{bfls}, \cite{RW}, \cite{M}, \cite{fls}). Mathematicians often consider \Z-graded,
rather than \zt-graded spaces. Since signs are determined only by the induced \zt-grading,
\Z-graded \linf\ algebras are examples of \zt-graded \linf\ algebras; however there are some
important differences in the classification, not only due to the fact that only some of the \zt-graded
\linf\ structures are \Z-graded \linf\ algebras, but also because the \Z-graded automorphisms are
a subgroup of the \zt-graded automorphisms, so the equivalence classes are potentially
quite different. Examples of \Z-graded \linf\ algebra structures were 
considered in \cite{D}.

In the physics literature, one usually considers $\Z_2$-graded spaces.
That's the case in our consideration also: throughout this paper, all
spaces will be $\Z_2$-graded, and we will work in the parity reversed
definition of the \linf\ structure. In \cite{fp2} we classified all
\linf\ algebras of dimension less than or equal to 2, in \cite{fp3}
constructed miniversal deformations for all \linf\ structures on a
space of three odd dimension - which
correspond to ordinary Lie algebras, and in \cite{fp4} we gave a
complete classification of all \linf\ algebras of dimension $1|2$.

The picture in the $2|1$ dimensional case is more complicated than for
$1|2$-dimensional algebras, because the space of $n$-cochains on a
$1|2$ dimensional space has dimension $6|6$ for $n>1$, while the space
of $n$-cochains on a $2|1$- dimensional space has dimension
$3n+2|3n+1$, making it more difficult to classify the nonequivalent
structures.  In this case we gave in \cite{bfp1} a complete
classification of only those \linf\ algebras which correspond to degree
1 and degree 2 coderivations.

In this paper, based on the classification in \cite{bfp1}, we classify
\linf\ algebras which are extensions
of degree 1 and degree 2 coderivations, in other words,  extensions of
 \zt-graded Lie algebras
as \linf\ algebras.  We introduce here all the
necessary results in order to make this article self contained.

\section{Basic Definitions}

\subsection{\linf\ algebras}

We work in the framework of the parity reversion $W=\Pi V$ of the usual
vector space $V$ on which an \linf\ algebra structure is defined,
because in the $W$ framework,  an \linf\ structure is simply an odd
coderivation $d$ of the symmetric coalgebra $S(W)$, satisfying
$d^2=0$,  in other words, it is an odd codifferential in the \zt-graded
Lie algebra of coderivations of $S(W)$.  As a consequence, when
studying \zt-graded Lie algebra structures on $V$, the parity is
reversed,  so that a $2|1$-dimensional vector space $W$ corresponds to
a $1|2$-dimensional \zt-graded Lie structure on $V$.  Moreover,  the
\zt-graded anti-symmetry of the Lie bracket on $V$ becomes the
\zt-graded symmetry of the associated coderivation $d$ on $S(W)$.

A formal power series $d=d_1+\cdots$, with $d_i\in L_i=\hom(S^i(W),W)$ determines
an element in $L=\hom(S(W),W)$, which is naturally identified with
$\coder(S(W))$, the space of coderivations of the symmetric
coalgebra $S(W)$. Thus $L$ is a \zt-graded Lie algebra.
An odd element $d$ in $L$ is called a  \emph{codifferential} if
$\br dd=0$. We also say that $d$ is
an \linf\ structure on $W$.A  detailed description of \linf\ algebras can be obtained in
\cite{lm,ls,pen3,pen4}.

If $g=g_1+\cdots\in\hom(S(W),W)$, and $g_1:W\ra W$ is invertible,
then $g$ determines a coalgebra automorphism of $S(W)$ in a natural
way, which we will denote by the same letter $g$. Moreover, every
coalgebra automorphism is determined in this manner. Two codifferentials
$d$ and $d'$ are said to be \emph{equivalent} if there is a coalgebra
automorphism $g$ such that $d'=g^*(d)=g\inv d g$.
An automorphism is said
to be linear when it is determined by a linear map $g_1$. Two codifferentials
are said to be \emph{linearly equivalent} when there is a linear equivalence
between them.
If $d$ and $d'$ are
codifferentials of a fixed degree $N$, then they are equivalent precisely
when they are linearly equivalent. Thus we can restrict ourself to linear automorphisms when
determining the equivalence classes of elements in $L_N$. In \cite{bfp1}, we classified
all codifferentials in $L_1$ and $L_2$. In this paper, our aim is to classify the extensions
of such codifferentials of a fixed degree to more general \linf\ algebras.

\subsection{Equivalent codifferentials and extensions}

We will use the following facts, which are
established in \cite{fp1}, to aid in the classification
of extensions of codifferentials of degree 1 and 2 to more
general \linf\ structures.

If $d$ is an \linf\ structure on $W$, and $d_N$ is the first
nonvanishing term in $d$, then $d_N$ is itself a codifferential, which we
call the \emph{leading term of $d$}, and we say that $d$ is
\emph{an extension} of $d_N$. Define the cohomology operator $D$
by $D(\ph)=\br\ph {d_N}$, for $\ph\in L$. Then the following
formula holds for any extension $d$ of $d_N$ as an \linf\
structure, and all $n\ge N$:
\begin{equation}\label{extension}
D(d_{n+1})=-\frac12\sum_{k=N+1}^{n}\br{d_k}{d_{n+N-k+1}}.
\end{equation}
Note that the terms on the right all have index less than $n+1$.
If a coderivation $d$ has been constructed up to terms of degree
$m$, satisfying \refeq{extension} for $n=1\dots m-1$, then the
right hand side of \refeq{extension} for $n=m$ is automatically a
cocycle. Thus $d$ can be extended to the next level precisely when
the cocycle given by the right hand side is trivial. There may be many
nonequivalent extensions, because the term $d_{n+1}$ which we add
to extend the coderivation is only determined up to a cocycle.  An
extension $d$ of $d_N$ is given by any coderivation whose leading term is
$d_N$, which satisfies
\refeq{extension} for every $m=N+1\dots$. The theory here is
parallel to the theory of formal deformations of an algebra structure;
the extension of a codifferential $d_N$ to a more complicated
codifferential $d$ resembles the process of extending an infinitesimal
deformation to a formal one.

Classifying the extensions of $d_N$ can be quite complicated.
However, the following theorem often makes it easy to classify the
extensions.
\begin{thm}\label{zerocohomology}
If the cohomology $H^n(d_N)=0$, for $n>N$, then any extension of
$d_N$ to a \linf\ structure $d$ is equivalent to the structure
$d_N$.
\end{thm}

Before classifying the extensions of a codifferential $d_N$, we
need to recall the classification of codifferentials in $L_N$ up to equivalence,
which was given in \cite{bfp1}. A \emph{linear automorphism} of $S(W)$ is an
automorphism determined by an isomorphism $g_1:W\ra W$. If $g$ is
an arbitrary automorphism, determined by maps $g_i:S^i(W)\ra W$,
and $W$ is finite dimensional, then $g_1$ is an isomorphism, so
this term alone induces an automorphism of $S(W)$ which we call
the \emph{linear part} of $g$.

The following theorem simplifies the classification of equivalence
classes of codifferentials in $L_N$.
\begin{thm}\label{simpleauto}
If $d$ and $d'$ are two codifferentials in $L_N$, and $g$ is an
equivalence between them, then the linear part of $g$ is also an
equivalence between them.
\end{thm}
Thus the equivalence classes of elements in $L_N$ are completely determined by
the action of the linear automorphisms on the coderivations.

We will also use the following result.
\begin{thm}\label{simpleequiv}
Suppose that $d$ and $d'$ are equivalent codifferentials.  Then
their leading terms have the same degree and are equivalent.
\end{thm}

As a consequence of these theorems, the classification of \linf\ structures
can be given as follows. First,  find all equivalence classes
of codifferentials of degree $N$. For each equivalence class, choose some representative
$d_N$ and determine the equivalence classes of extensions of the codifferential $d_N$.
The first part of this task was completed in \cite{bfp1}. Here we address the second
part.

Let us first establish some basic notation for the cochains.
Suppose $W=\langle w_1,w_2,w_3\rangle$, with $w_1$ an odd element and $w_2,w_3$
 even elements.
If $I=(i_1,i_2,i_3)$ is a
multi-index, with $i_1$ and $i_2$ either zero or one, let
$w_I=w_1^{i_1}w_2^{i_2}w_3^{i_3}$. For simplicity, we will sometimes abbreviate
$w_I$ by $I$.
Then for $n\ge 1$,
\begin{align*}
(S^{n}(W))_e=&\langle (0,p,n-p)|0\le p\le n\rangle,\qquad |(S^{n}(W))_e|=n+1\\
(S^{n}(W))_o=&\langle (1,q,n-q-1)|0\le q\le n-1\rangle,\qquad |(S^{n}(W))_0|=n
\end{align*}

If $\lambda$ is a
linear automorphism of $S(W)$, then in terms of the standard basis of $W$, its restriction to $W$ has matrix
\begin{equation}\label{linauto}
\lambda=\begin{pmatrix}
q&0&0\\
0&r&t\\
0&s&u\\
\end{pmatrix}
\end{equation}
where $q(ru-st)\ne 0$.
We will sometimes express $\lambda$ by the submatrix $\bigl(\begin{smallmatrix}r&t\\s&u\end{smallmatrix}\bigr)$.
It is useful to note that for a linear automorphism
$$\lambda(w_I)=\lambda(w_1)^{i_1}\lambda(w_2)^{i_2}\lambda(w_3)^{i_3},$$
so that
\begin{equation}\label{lambdaxy}
\lambda(z,x,y)=\sum_{i=0}^x\sum_{j=0}^y(z,i+j,x+y-i-j)\binom xi \binom yj q^zr^is^{x-i}t^ju^{y-j}.
\end{equation}

For a multi-index $I$, let $|I|=i_1+i_2+i_3$. For $n=|I|$, define $\ph^I\in L_n=\hom(S^n(W),W)$ by
$
\ph^I_j(w_J)=I!\delta^I_J w_j,
$
where $I!=i_1!i_2!i_3!$. Then
$L_n=\langle \ph^I_j, |I|=n \rangle$. If $\ph$ is odd, we will denote
it by the symbol $\psi$ to make it easier to distinguish the even
and odd elements. Then
\begin{align*}
(L_{n})_e&=\langle\ph^{1,q,n-q-1}_1,\ph^{0,p,n-p}_2,\ph^{0,p,n-p}_3|1\le q\le n-1,1\le p\le n\rangle\\
(L_{n})_o&=\langle\psi^{1,q,n-q-1}_2,\psi^{1,q,n-q-1}_3,\psi^{0,p,n-p}_1|1\le q\le n-1,1\le p\le n\rangle,
\end{align*}
so that $|L_n|=3n+2|3n+1$.

A linear automorphism $\lambda$ acts on the space of coderivations by
\begin{equation*}\lambda^*(\ph)=\lambda\inv\ph\lambda.\end{equation*}
A somewhat complicated formula for the action can be given as follows.
\begin{align*}
\lambda^*(\ph^{z,a,b}_1)=&\sum_{\substack{x=0\cdots a+b\\i=x-b\cdots a}}
\phd{z,x,a+b-x}1q\inv V_i \\
\lambda^*(\ph^{z,a,b}_2)=&\sum_{\substack{x=0\cdots a+b\\i=x-b\cdots a}}
(\phd{z,x,a+b-x}2\tfrac u\Delta-\phd{z,x,a+b-x}3\tfrac s\Delta)V_i\\
\lambda^*(\ph^{z,a,b}_3)=&\sum_{\substack{x=0\cdots a+b\\i=x-b\cdots a}}
(-\phd{z,x,a+b-x}2\tfrac t\Delta+\phd{z,x,a+b-x}3\tfrac r\Delta)V_i,
\end{align*}
where $V_i=r^is^{x-i}t^{a-i}u^{b+i-x}\binom x i \binom {a+b-x}{a-i}$ and $\Delta=ru-st$.

\section{Classifying Codifferentials}


In \cite{bfp1}, we showed that codifferentials of degree $n$ fall into two distinct
families,  those of the
\emph{first kind}
\begin{equation}
\sum_{q=0}^{n-1}\psi^{1,q,n-q}_2b_q+\psi^{1,q,n-q}_3c_q.
\end{equation}
and those of the \emph{second kind}
\begin{equation}
d=\sum_{p=0}^n\psi^{0,p,n-p}_1a_p.
\end{equation}

Moreover, any expression of either kind gives a codifferential.  Thus we have determined all codifferentials
of degree $N$.  However, the process of classification requires that we determine the equivalence classes of
codifferentials under the action of the automorphism group of the symmetric coalgebra, and we are a long
way away from this classification at this stage.

Let us call the degree of the leading term of a codifferential the
\emph{order} of that codifferential.
In \cite{bfp1} we classified all codifferentials of degree 1 and 2. Here
we will recall the relevant details, as a starting point to study
the extensions.

\section{Codifferentials of Degree 1 and Their Extensions}

 Let us
suppose that $d$ is an odd, degree 1 codifferential of the first kind.
Then it is equivalent to $\psi^{1,0,0}_2$. Its cohomology is as
follows.  We define the coboundary operator $D$ by $D(\ph)=\br{\ph}d$.
Then computing brackets, we see that \begin{align*}
D(\ph^{0,p,n-p}_2)&=\psi^{1,p-1,n-p}_2p,\quad
D(\psi^{0,p,n-p}_1)=\ph^{1,p-1,n-p}_1p+\ph^{0,p,n-p}_2\\
D(\ph^{0,p,n-p}_3)&=\psi^{1,p-1,n-p}_3p,\quad D(\psi^{1,q,n-q-1}_2)=0\\
D(\ph^{1,q,n-q-1}_1)&=-\psi^{1,q,n-q-1}_2,\quad D(\psi^{1,q,n-q-1}_3)=0
\end{align*}
From this, it follows easily that the cohomology of $d$ is zero, so by
\thmref{zerocohomology} all extensions of $d$ are equivalent to
$d$. This completes the picture for codifferentials of the first kind
of degree 1.

If $d$ is a codifferential of the second kind of degree 1, it is equivalent
 to $d'=\psi^{0,1,0}_1$.
Calculating coboundaries, we have
\begin{align*}
D(\ph^{0,p,n-p}_2)=&-\psi^{0,p,n-p}_1,&
D(\psi^{0,p,n-p}_1)=&0\\
D(\ph^{0,p,n-p}_3)=&0,&
D(\psi^{1,q,n-q-1}_2)=&\ph^{1,q,n-q-1}_1+\ph^{0,q+1,n-q-1}_2\\
D(\ph^{1,q,n-q-1}_1)=&\psi^{1,q+1,n-q-1}_1,&
D(\psi^{1,q,n-q-1}_3)=&\ph^{0,q+1,n-q-1}_3
\end{align*}

The cohomology of this codifferential is also equal to zero.
Thus every extension of a codifferential
of the second kind is equivalent to the original codifferential.
This completes the classification of all codifferentials whose leading term has degree 1. They are
all equivalent to the codifferential given by the leading term alone.  Thus there are no interesting
extensions of degree 1 codifferentials.

\section{Codifferentials of Degree 2}
We now consider how to extend a codifferential $d_2$ of degree 2 to a more general codifferential.
For a codifferential $d$, let $\aut(d)$ be the subgroup of automorphisms of $S(W)$ fixing $d$,
and $\laut(d)$ be the subgroup of $\aut(d)$
consisting of linear automorphisms. If $d=d_1\mplus d_k$ is a codifferential consisting of terms only up
to degree $k$, then a generalized automorphism $g$ of $d$ is an automorphism of $S(W)$ such that
$g^*(d)=d+\ho$. The set $\gaut(d)$ of generalized automorphisms of $d$ is a subgroup of the automorphisms
of $S(W)$ containing $\aut(d)$.

Let $d_2$ be a codifferential of degree $2$.
The groups $\aut(d_2)$ and $\laut(d_2)$ act on the set of nonzero cohomology classes of $d_2$. We shall say
that two cohomology classes $\delta_1$ and $\delta_2$ are equivalent if there is an element $f$ in $\aut(d_2)$ such
that $f^*(\delta_1)=\delta_2$,  and linearly equivalent if $f$ lies in $\laut(d_2)$.

Any automorphism $g$ can be expressed in the form $$g=\lambda\prod_{k=2}^\infty\exp(\alpha_k),$$ where
$\lambda$ is linear $\alpha_k\in L_k$ is a coderivation.  We call $\lambda$ the linear part of $g$. An
automorphism of the form $g=\prod_{k=2}^\infty\exp(\alpha_k)$, is called a formal automorphism of $S(W)$.
The set $\autf(S(W))$ offormal automorphisms is a normal subgroup in $\aut(S(W))$, so the linear part $\lambda$
of a automorphism $g$ is uniquely determined.

The action of $\aut(S(W))$ on the space $L$ of coderivations of $S(W))$ is given by
$g^*=\prod_{k=\infty}^2\exp(-\ad_{\alpha_k})\lambda^*$, where
$\ad_{\alpha_k}(\beta)=[\alpha_k,\beta]$, and $\lambda^*(\beta)=\lambda\inv\beta\lambda$.

Suppose that $d=d_2 + d_{k}+ d_{k+1}+\cdots$ is an extension of $d_2$.  The following theorem will help
to classify such extensions.
\begin{thm}\label{equivclasses}
Suppose that
$d=d_2 + d_{k}+ d_{k+1}+\cdots$ is a codifferential. Then $d_{k}$ is a cocycle with respect to
the coboundary operator $D=\br{\bullet}{d_2}$. Moreover,
\begin{enumerate}
\item $d$ is equivalent to a codifferential whose first nonzero term
after $d_2$ is of higher order than $k$ iff $d_{k}$ is a $D$-coboundary,.
\item If $d'=d_2+ d'_{k}+ d'_{k+1}+\cdots$ is equivalent to $d$, then the cohomology classes of
$d'_{k}$ and $d_k$ are  linearly
equivalent.
\end{enumerate}
\end{thm}
\begin{proof}
Since $[d,d]=2D(d_k)+\cdots$, it follows that $D(d_k)=0$. If $d_k=D(\alpha)$, then applying $g=\exp(\alpha)$
to $d$, we obtain $g*(d)=d_2+d_k-D(\alpha)+\cdots$, so we have eliminated the term $d_k$. On the other hand,
suppose that $g^*(d)=d_2+ d'_{k}+ d'_{k+1}+\cdots$, for some $g=\lambda\exp(\beta_2)\cdots$.
Since we must have $\lambda^*(d_2)=d_2$, we compute
\begin{align*}
g^*(d)=&d_2+\lambda^*(d_k)-D(\beta_2)-\cdots-D(\beta_{k-1})
\\&+\tfrac12[\beta_2,D(\beta_2)]+\cdots
\end{align*}
It follows that $D(\beta_i)=0$ for $i<k-1$ and $-D(\beta_{k-1})+\lambda^*(d_k)=d'_k$.
If $d_k'=0$, this says that $\lambda^*(d_k)$ is a coboundary, and since the coboundary map commutes with
automorphisms, we see that $d_k$ is also a coboundary. It is also clear that the cohomology class of
$\lambda^*(d_k)$ coincides with that of $d_k'$.
\end{proof}
Let us say that a codifferential $d$ is standard if it is of the form $d_2 + d_{k}+ d_{k+1}+\cdots$,
where $d_{k}$ is a nontrivial cocycle for $d_2$.  By the first part of this theorem,  every nontrivial
extension of $d_2$ is equivalent to a standard codifferential.  For a codifferential in standard form,
let us refer to the cohomology class of $d_{k}$ as the secondary term of $d$.

In general, we don't expect $d_2 +d_{k}$ to be a codifferential for an extension
$d=d_2 + d_{k}+ d_{k+1}+\cdots$ of $d_2$.  However,
for the examples  which arise in this paper, it turns out to be true.
We state a theorem which is useful in characterizing the extensions of a codifferential of the
form $d=d_2 +d_{k}$.
\begin{thm}\label{equivext}
Suppose $d=d_2 +d_{k}$ is a codifferential. Let $D_2$ be the coboundary operator determined by $d_2$,
and $D_k$ be the one given by $d_{k}$.
Let $d'=d_2 +d_{k}+d_{l}+\cdots$ be an extension of $d$, where $l>k$. 
Then
\begin{enumerate}
\item
$d_l$ is a $D_2$-cocycle. Moreover $D_k(d_l)$ is a $D_2$-coboundary.
\item
If $d_l$ is a $D_2$-coboundary, then $d'$ is equivalent to an extension  whose third term
 has degree larger
than $l$.
\item If $H^n(d_2)=0$ for $n>k$, then any extension $d'$ of $d$ is equivalent to $d$.
\item
If $d_l=D_k(\eta)$ for some $D_2$-cocycle $\gamma$,
then $d'$ is equivalent to an extension whose third term has degree larger
than $l$.
\item If every $D_2$-cocycle $\delta$ whose order is larger than $k$
for which $D_k(\delta)$ is a $D_2$-coboundary is of the form
$\delta=D_k(\eta)$ for some $D_2$-cocycle $\eta$, then any extension of $d$
is equivalent to $d$.
\end{enumerate}
\begin{proof}
The lowest order term in $\tfrac12[d',d']$ is $D_2(d_l)$, and therefore it must be equal to zero.
In fact, $D_2(d_l)=\cdots=D_2(d_{k+l-1})=0$, and $D_2(d_{k+l})=-D_k(d_l)$, by matching the order of terms,
so the second assertion holds. If $d_l=D_2(\gamma)$, applying $\exp(-\ad_\gamma)$ to $d'$ removes the $d_l$
term. Similarly, if $d_l=D_k(\gamma)$ for a $D_2$ cocycle $\gamma$, then applying $\exp(-\ad_\gamma)$ to
$d_l$ removes the $d_l$ term. Since $d_l$ satisfies the property that its $D_k$-coboundary is a $D_2$ coboundary,
the hypotheses of part 5 are met, so that $d_l$ satisfies the requirements of part 4, and it can be eliminated.
Thus all higher order terms can be successively eliminated, and $d'$ is equivalent to $d$.
\end{proof}
\end{thm}
A useful generalization of this theorem is as follows.
\begin{thm}\label{corollary} Let $d_e=d_2+\delta$ be a codifferential,
where $\delta=d_3\mplus d_k$ and $D_2$, $D_e$ be the coboundary operators
determined by $d_2$ and $d_e$, respectively. Suppose that every odd $D_2$-cocycle $\alpha$ of degree
greater than $k$ extends to a $D_e$-cocycle, and that $H^n_o(d_e)=0$ for $n>k$. Then every extension
of $d_e$ is equivalent to $d_e$.
\end{thm}
\begin{proof}
Let $d'=d_e+d_l+\ho$ be an extension of $d_e$. By the previous theorem, we know that $d_l$ is a $D_2$-cocycle,
so it extends to some $D_e$-cocycle $\delta$. Since the cohomology of $D_e$ vanishes for $n>k$,
$\delta=D_e(\gamma)$ for some cochain $\gamma$.  Now $d'=d_e+\delta+\ho$, and thus applying $\exp(-\ad_\gamma)$ to
$d'$ will eliminate the term of degree $l$ in $d'$.
\end{proof}
For each of the equivalence classes of degree 2 codifferentials,
we will study the nontrivial extensions.  Such extensions exist only
when the odd part of $H^n$ does not vanish for some $n>2$, because the hypotheses of \thmref{corollary} above are met
for trivial reasons, because every codifferential is a sum of odd cochains of the same kind, and the brackets
of cochains of the same type with each other always vanish, so every $D_2$-cocycle is automatically a $D_e$-cocycle.
Ordinarily, when considering extensions,  we have to construct them term by term, because
a finite number of terms may not determine a codifferential.  However, in our examples, because the cocycles
all have trivial brackets with respect to each other, so there is never any question about whether
adding a cocycle gives a codifferential.

Let us state a conjecture which we will use to classify the codifferentials.
Two codifferentials are said to be formally equivalent if there is a formal automorphism expressing an
equivalence between them.
\begin{con}\label{conjecture}
Suppose that $d=d_2+d_k$ is a codifferential, $d_e=d_2+d_k+d_l+\cdots$ and  $d_e'=d_2+d_k+d_l'+\cdots$.
If $d_e'$ is formally equivalent to $d_e$, then $d'_l-d_l$ is the leading term of a $D$-coboundary with respect to
the coboundary operator $D=[\bullet, d]$.
\end{con}
In particular,  this conjecture implies that if $d_l$ is not the leading term of a coboundary with
respect to $D$, then $d_e$ is not formally equivalent to $d$. Of course, we are actually interested
in whether two codifferentials are equivalent, not formally equivalent.  However, we can often use the
following fact to reduce the study of equivalent extensions of a codifferential $d$ to the study
of the action of $\laut(d)$ on the formal equivalence classes of extensions of $d$.

Suppose that $d_k$ is a $d_2$-cocycle. If every linear automorphism of $d_2$ which takes $d_k$
to a $d_2$-cohomologous cochain fixes $d_k$, then the linear part of such an automorphism fixes $d_k$.
In other words, the linear part of a generalized automorphism of $d=d_2+d_e$ is a generalized automorphism
of $d$.  This statement is clear because of \thmref{equivclasses}.

Another useful observation is the following. A diagonal linear automorphism acts on cochains
by mapping them to multiples of themselves.
Since a multiple of a nontrivial cocycle cannot be cohomologous to the cocycle unless it is equal to the
cocycle, if $\lambda$ is the linear part of a generalized automorphism of $d=d_2+d_e$, then $\lambda$
fixes $d_e$.

We will start with studying codifferentials of degree 2 of the second
kind.

\section{Codifferentials of Degree 2 of the Second Kind and Their
Extensions}

There are only two types of codifferentials up to equivalence:
$\psd{0,2,0}1$ and $\psd{0,1,1}1$. For $\psd{0,1,1}1$, it was shown 
(see \cite{bfp1}) that  the odd cohomology $H^n_o$ vanishes if $n>1$.
Thus there are no notrivial extensions.  For $d=\psd{0,2,0}1$, the
situation is a bit more complicated. Let
$D(\ph)=\br{\ph}{\psd{0,2,0}1}$. Then we obtain the following table of
coboundaries.
\begin{align*} D(\phd{1, q,
n-q-1}1) &=\psd{0, 2+q, n-q-1}1\\
D(\phd{0, p, n-p}2) &=-2\psd{0,p+1,n-p}1\\
D(\phd{0, p, n-p}3) &=0\\
D(\psd{0, p, n-p}1) &=0\\
D(\psd{1,q, n-q-1}2) &=2\phd{1,q+1,n-q-1}1+\phd{0,q+2,n-q-1}2\\
D(\psd{1, q,n-q-1}3) &=\phd{0,q+2,n-q-1}3\\
\end{align*}
It is not hard to see that the cohomology is given by
\begin{align*}
&H^1=\langle \psd{0,0,1}1,\psd{0,1,0}1,\phd{0,0,1}3,\phd{0,1,0}3,2\phd{1,0,0}1+\phd{0,1,0}2\rangle\\
&H^n=\langle \psd{0,0,n}1,\phd{0,0,n}3,\phd{0,1,n-1}3,2\phd{1,0,n-1}1+\phd{0,1,n-1}2,\rangle,\quad\text{if $n>1$}
\end{align*}

Let us now study extensions of $d$.
Suppose that $\lambda$ is a linear automorphism of $d$.  It is easy to check that the condition
$\lambda^*(d)=d$ is equivalent to $t=0$ and $r^2=1$ in the standard expression (\refeq{linauto}) for $\lambda$.

For $k>2$, we can extend $d$  to the codifferential
\begin{equation}
d_e=\psd{0,2,0}1+\psd{0,0,k}1a.
\end{equation}
When $a$ is a nonzero number, by using a diagonal automorphism of $d$,
we see that $d_e$ is equivalent to the codifferential $\psd{0,2,0}1+\psd{0,0,k}1$, so we may assume
that $a=1$. If $\lambda$ is an automorphism of $d$, then
\begin{equation}
\lambda^*(\psd{0,0,k}1)=\sum_{x=0}^k\psd{0,x,k-x}1\frac{s^xu^{k-x}}q.
\end{equation}
In order for $\lambda$ to be the linear part of a generalized automorphism of $d_e$, $\lambda^*(\psd{0,0,k}2)$
must be $D$-cohomologous to $\psd{0,0,k}2$.  This can only happen if $u^k/q=1$. Thus, in this case, there
are examples of generalized automorphisms of $d_e$ for which the linear part of the automorphism does not
fix $d_e$.

We compute
\begin{align*}
&D_e(\phd{0,0,n}3)=-\psd{0,0,k+n-1}1k\\
&D_e(\phd{0,1,n-1}3)=-\psd{0,1,k+n-2}1k\\
&D_e(2\phd{1,0,n-1}1+\phd{0,1,n-1}2)=2\psd{0,0,k+n-1}1
\end{align*}
The first and the third even $D$-cocycle combine to give the cocycle
$2\phd{0,0,n}3+2\phd{1,0,n-1}1k+\phd{0,1,n-1}2k$, while the second extends to the $D_e$-cocycle
$2\phd{0,1,n-1}3-\phd{0,0,n+k-2}2$.
Also, we note that $\psd{0,0,n}1$ is a coboundary for $n\ge k$.

The cohomology of $d_e$ is given by
\begin{align*}
H^1&=\langle \psd{0,0,1}1,\psd{0,1,0}1,2\phd{0,0,1}3+2\phd{1,0,0}1k+\phd{0,1,0}2k,
2\phd{0,1,0}3-\phd{0,0,k-1}2\rangle\\
H^n&=\langle \psd{0,0,n}1,2\phd{0,0,n}3+2\phd{1,0,n-1}1k+\phd{0,1,n-1}2k,2\phd{0,1,n-1}3-\phd{0,0,n+k-2}2\rangle,\\
&\quad\text{if $1<n<k$}\\
H^n&=\langle 2\phd{0,0,n}3+2\phd{1,0,n-1}1k+\phd{0,1,n-1}2k,2\phd{0,1,n-1}3-\phd{0,0,n+k-2}2\rangle,\\
&\quad\text{if $n\ge k$}
\end{align*}
Thus we have no higher order odd cohomology, so every extension of $d_e$ is equivalent to it, by
\thmref{corollary}.  In this case, we did not need to study the action of the automorphism group
on $d_e$ in order to classify the extensions. But in some later examples we will need to study
this action very carefully.

\section{Codifferentials of Degree 2 of the First Kind and Their
Extensions}

Let us label the codifferentials representing the equivalence classes of degree 2 codifferentials of
the first kind as follows:(see \cite{bfp1} for a proof that these give a complete set of representatives
of the equivalence classes.)
\begin{align*}
&\dsharp=\psd{1,1,0}2+\psd{1,1,0}3+\psd{1,0,1}3\\
&\dt{c}=\psd{1,1,0}2+\psd{1,0,1}3c, |c|\le 1\\
&\dstar=\psd{1,1,0}3
\end{align*}

The cohomology of $\dsharp$ vanishes if $n>2$, so there are no nontrivial extensions. The cohomology
of $\dt c$ depends on the value of $c$, so there are several cases to consider.

\subsection{The Codifferential $\dt c=\psd{1,1,0}2+\psd{1,0,1}3c$}
For an arbitrary value of $c$, the coboundaries are given by
\begin{align*}
D(\phd{1, q, n-q-1}1) &=-\psd{1, 1+q, n-q-1}2-\psd{1, q, n-q}3c\\
D(\phd{0, p, n-p}2) &=\psd{1,p, n-p}2(p-1+c(n-p))\\
D(\phd{0, p, n-p}3) &=\psd{1,p, n-p}3(p+c(n-p-1))\\
D(\psd{0, p, n-p}1) &=\phd{1,p,n-p}1(p+c(n-p))+\phd{0, 1+p, n-p}2+\phd{0, p, n-p+1}3c\\
D(\psd{1, q, n-q-1}2) &=0\\
D(\psd{1, q, n-q-1}3) &=0
\end{align*}

First consider what the linear automorphism group of $\dt c$ consists of.
If $\lambda$ is a linear automorphism of $S(W)$ is
given by the standard form (\refeq{linauto}) then
\begin{align*}
\lambda^*(d_c)=
&\psd{1,1,0}2\pfrac{q(ru-cst)}{ru-st}+\psd{1,1,0}3\pfrac{(c-1)qrs}{ru-st}\\
+&\psd{1,0,1}2\pfrac{(1-c)qtu}{ru-st}+\psd{1,0,1}3\pfrac{q(cur-st)}{ru-st}.
\end{align*}
If $c=1$ there are no nontrivial extensions, because the higher cohomology vanishes.
Assuming that $c\ne1$,
we observe that $\lambda^*(\dt c)=\dt c$ in only two cases.
First, we could have $s=t=0$ and $q=1$.  Then $\lambda$ is diagonal.
Secondly,  we could have $u=r=0$, $q=-1$ and $c=-1$.  Thus the second case only occurs for the special value $c=-1$.
This special case is of particular interest to us, so we will consider it separately.

For generic values of $c$, the higher cohomology vanishes, so we don't have any nontrivial extensions.
Let us consider the special cases for which the higher cohomology does not vanishes.
\subsubsection{Extensions of $\dt c$, when $1/c>1$ is a positive integer}
Let $m=1/c$. Then we have
\begin{align*}
&H^1=\langle \psd{1,0,0}2,\psd{1,0,0}3,\phd{0,1,0}2,\phd{0,0,1}3\rangle\\
&H^2=\langle \psd{1,0,1}3\rangle\\
&H^m=\langle \phd{0,0,m}2\rangle\\
&H^{m+1}=\langle \psd{1,0,m}2\rangle\\
&H^n=0,\qquad\text{otherwise}
\end{align*}
except when $c=1/2$, in which case, since $m=2$, $H^2=\langle \psd{1,0,1}3,\phd{0,0,2}2\rangle$.
Of course, the case $c=1/2$ is not interesting because there is no higher cohomology, so let us
assume $c\ne 1/2$.

A nontrivial extension of
$\dt c$ must be of the form
\begin{equation*}
d_e=\psd{1,1,0}2+\psd{1,0,1}3c+\psd{1,0,m}2a
\end{equation*}

It is easy to see, using a diagonal linear transformation, that we can take $a=1$.
Note that a linear automorphism $\lambda$ of $d$ preserves $d_e$ precisely when
$u^m=r$. (Remember that $\lambda$ is diagonal and $q=1$.)

Let us study the cohomology of $d_e$. First,  note that
\begin{align*}
D_e(\phd{0,1,0}2)=&\psd{1,0,m}2\\
D_e(\phd{0,0,1}3)=&-\psd{1,0,m}2m,
\end{align*}
so that the two 1-cohomology
classes for $\dt c$ are replaced by the single $D_e$-cohomology class $\phd{0,1,0}2m+\phd{0,0,1}3$, and
$\psd{1,0,m}2$ becomes a $D_e$-coboundary.  Note that $\phd{0,0,m}2$ remains a $D_e$-cocycle.
Thus the cohomology of $D_e$ is given by
\begin{align*}
&H^1=\langle \psd{1,0,0}2,\psd{1,0,0}3,\phd{0,1,0}2m+\phd{0,0,1}3\rangle\\
&H^2=\langle \psd{1,0,1}3\rangle\\
&H^m=\langle \phd{0,0,m}2\rangle\\
&H^n=0,\qquad\text{otherwise}
\end{align*}
By \thmref{corollary}, it follows that any extension of $d_e$ is equivalent to $d_e$, so we have
found all nonequivalent extensions of $\dt c$. Notice that in this case as well, the automorphism
group of $d_e$ plays no role in the classification.

\subsubsection{Extensions of $\dt 0$}
The cohomology is given by
\begin{align*}
&H^1=\langle \psd{1,0,0}2,\psd{1,0,0}3,\phd{0,1,0}2,\phd{0,0,1}3\rangle\\
&H^n=\langle \psd{1,0,n-1}3,\phd{0,0,n}3\rangle,\qquad\text{if $n>1$}
\end{align*}

If we take an extension of the form $\psd{1,1,0}2+\psd{1,0,k}3a$, then as usual, we can assume $a=1$,
so we may as well assume that our extended codifferential is
$$
d_e=\psd{1,1,0}2+\psd{1,0,k}3.
$$
The linear part $\lambda$ of a generalized automorphism of $d_e$ is a diagonal automorphism with
$q=1$ and $u^{k-1}=1$, so that in fact $\lambda^*(d_e)=d_e$.
We note that
\begin{align*}
&D_e(\phd{0,1,0}2)=0\\
&D_e(\phd{0,0,n}3)=\psd{1,0,k+n-1}3(n-k)
\end{align*}
so that $\phd{0,1,0}2$ remains a $D_e$-cocycle, but with the exception of $n=k$, the cochains $\phd{0,0,n}3$
give rise to $D_e$-coboundaries, and are no longer cocycles.  Also, $\psd{1,0,m}3$ is a $D_e$-coboundary
for $m\ge k$ except for $m=2k-1$. Thus, after applying an appropriate automorphism,
we can assume that the third order term in any nontrivial extension of $d_e$ is of the form $\psd{1,0,2k-1}3a$,
for some nonzero $a$, and we have the following candidate for a possible nontrivial extension of $d_e$.
$$
\dzeroa=\psd{1,1,0}2+\psd{1,0,k}3+\psd{1,0,2k-1}3a.
$$
Note that for any generalized automorphism of $\dzeroa$, its linear part $\lambda$ fixes $d_e$, so that it is a diagonal
automorphism with $u^{k-1}=1$. But then it follows that $\lambda^*(\psd{1,0,2k-1}3a)=\psd{1,0,2k-1}3a$,
so that $\lambda^*$ fixes $\dzeroa$.
If $\dzeroa$ were equivalent to $d_e$, $\psd{1,0,2k-1}3$
would be the leading term of a coboundary with respect to $D_e$, by Conjecture 5.4, and since this is not
true, we see that $\dzeroa$ is a nontrivial extension of $d_e$.  Moreover, two codifferentials with different
values of $a$ are not equivalent.

Because of this, we would expect
that $\psd{1,0,2k-1}3$ remains a cohomology class for $D_{0,a}=\br{\bullet}{\dzeroa}$, and this fact is
easily checked.  Moreover, $\phd{0,1,0}2$ still remains a $D_{0,a}$-cocycle.  Thus the cohomology for $\dzeroa$
is given by
\begin{align*}
&H^1=\langle \psd{1,0,0}2,\psd{1,0,0}3,\phd{0,1,0}2\rangle\\
&H^n=\langle \psd{1,0,n-1}3\rangle,\qquad\text{if $1<n\le k$}\\
&H^{2k}=\langle \psd{1,0,2k-1}3\rangle\\
&H^n=0,\qquad\text{otherwise}
\end{align*}
By \thmref{corollary}, we again see that any extension of $\dzeroa$ is equivalent to $\dzeroa$. Thus
we have classified all extensions of $\dt 0$.

\subsubsection{Extensions of $\dt c$, when $c\ne -1$ is a negative rational number}
We may assume that $-1<c<0$.
If we express $\frac c{c-1}=\frac rs$, as a fraction in lowest terms,
then $0<\frac rs<1/2$
Then the cohomology is given by the following table.
\begin{align*}
&H^1=\langle \psd{1,0,0}2,\psd{1,0,0}3,\phd{0,1,0}2,\phd{0,0,1}3\rangle\\
&H^2=\langle \psd{1,0,1}3\rangle\\
&H^{ks+1}=\langle \phd{0,kr,k(s-r)+1}3\rangle\\
&H^{ks+2}=\langle \psd{1,kr,k(s-r)+1}3\rangle\\
&H^n=0,\qquad\text{otherwise}
\end{align*}
A nontrivial odd cohomology class can be represented by a cocycle of
the form.
$\psd{1,kr,k(s-r)+1}3a$. If we add such a cocycle to $\dt c$, then by applying a linear automorphism,
one sees that up to equivalence, we may assume that $a=1$. Thus we have a nontrivial extension
\begin{equation}
d_e=\psd{1,1,0}2+\psd{1,0,1}3c+\psd{1,kr,k(s-r)+1}3.
\end{equation}
We obtain
\begin{equation}
D_e(\phd{0,lr,l(s-r)+1}3)=\psd{1,(l+k)r,(l+k)(s-r)+1}3(l-k)(s-r),
\end{equation}
which does not vanish unless  $l=k$, so we get the following table for the cohomology of $d_e$.

\begin{align*}
&H^1=\langle \psd{1,0,0}2,\psd{1,0,0}3,\phd{0,1,0}2(r-s)+\phd{0,0,1}3r\rangle\\
&H^{ms+2}=\langle \psd{1,mr,m(s-r)+1}3\rangle\qquad\text{if $0\le m<k$}\\
&H^{ks+1}=\langle \phd{0,kr,k(s-r)+1}3\rangle\\
&H^{2ks+2}=\langle \psd{1,2kr,2k(s-r)+1}3\rangle\\
&H^n=0,\qquad\text{otherwise}
\end{align*}
Thus the only candidate for a nontrivial extension of $d_e$ is
\begin{equation}
\dt {c,a}=\psd{1,1,0}2+\psd{1,0,1}3c+\psd{1,kr,r(s-r)+1}3+\psd{1,2kr,2k(s-r)+1}3a,
\end{equation}
where $a$ is a parameter which, as we will see, cannot be eliminated.
First,  note that $\phd{0,1,0}2(r-s)+\phd{0,0,1}3r$ remains a cocycle for $\dt{c,a}$.
Let $\ph=\phd{0,kr,k(s-r)+1}3$.
Since $D_{c,a}(\ph)=-\psd{1,3kr,3k(s-r)+1}3k(s-r)a$, $\ph$ no longer
generates $H^{ks+1}$.  However, $\ph$ extends to a cocycle $\ph'$, which is given by a power series
with leading term $\ph$.  Recall that for a codifferential which is not of fixed degree, the
spaces $H^n$ do not make sense in the usual manner; rather $H$ inherits a filtration from the natural
filtration on $L$, and $H^n$ is a quotient space of the $n$-th filtered part by the $n+1$-st, and in
this sense,  we have $H^{ks+1}=\langle \ph'\rangle$.  Note that other than this change,
all of the cohomology for $\dt{c,a}$ remains the same as for $d_e$.

If $\lambda^*(\dt c)=\dt c$, we have $\lambda=\diag(1,x,u)$, and we
can compute that
$$\lambda^*(\psd{1,kr,k(s-r)+1}3)=\psd{1,kr,k(s-r)+1}3x^{kr}u^{k(s-r)}.$$ 
This cochain is cohomologous to $\psd{1,kr,k(s-r)+1}3$ precisely when
it is equal to it, so that $x^{kr}u^{k(s-r)}=1$. But then it is
immediate that $\lambda^*(\dt{c,a})=d_{c,a}$. Thus we may conclude, 
from our conjecture, that $\dt {c,a}$ is a nontrivial extension of
$\dt c$ for each $a$, and the extensions are not equivalent for
different values of $a$.  Since the cohomology of $\dt c$ vanishes in
degree higher than $ks+1$, it follows that we have classified all
nonequivalent extensions of $\dt c$.

\subsubsection{Extensions of $d_{-1}$}
In this case $r=1$ and $s=2$, which makes things a bit more symmetric.
It turns out to be convenient to use different representatives for the
cohomology classes. Instead of using $\psd{1,k,k+1}3$ as the
representative for the basis of $H^{2k+2}$, if we choose
$\alpha_k=\psd{1,k+1,k}2+\psd{1,k,k+1}3$, then it is easily seen that
this element is a cocycle which is not a coboundary, and thus it can
be used as a basis element. Similarly, we can choose
$\beta_k=\phd{0,k+1,k}2+\phd{0,k,k+1}3$ as a basis for $H^{2k+1}$.

There are two types of linear automorphisms which preserve $d$. If
$\lambda=\diag(1,r,u)$, then $\lambda^*(\alpha_k)=\alpha_k(ru)^k$. On
the other hand, when
$\lambda=\Bigl(\begin{smallmatrix}-1&0&0\\0&0&t\\0&s&0\end{smallmatrix}\Bigr)$,
then $\lambda^*(\alpha_k)=-\alpha_k(st)^k$. Thus in both cases,
$\lambda^*(\alpha_k)$ is just a multiple of $\alpha_k$. Any nontrivial
extension of  $\dt{-1}$ is equivalent to one of the form
\begin{equation*}
d_e=\psd{1,1,0}2-\psd{1,0,1}3+\alpha_k.
\end{equation*}
Then
\begin{equation*}
D_e(\beta_l)=\alpha_{k+l}2(l-k),
\end{equation*}
which does not vanish unless $k=l$.  Thus, like the other cases when $c$ is a negative rational number, we have
a simple expression for the cohomology of $d_e$.  The table for the cohomology of $d_e$ is given by
\begin{align*}
&H^1=\langle \psd{1,0,0}2,\psd{1,0,0}3,-\phd{0,1,0}2+\phd{0,0,1}3\rangle\\
&H^{2(m+1)}=\langle\alpha_m\rangle\qquad\text{if $0\le m<k$}\\
&H^{2k+1}=\langle \beta_k\rangle\\
&H^{2(2k+1)}=\langle\alpha_{2k}\rangle\\
&H^n=0,\qquad\text{otherwise}
\end{align*}
Thus we have a single candidate for a nontrivial extension of $d_e$, given by
\begin{equation*}
d_{e,a}=\psd{1,1,0}2-\psd{1,0,1}3+\alpha_k+\alpha_{2k}a.
\end{equation*}
Now, we have already seen that any linear automorphism $\lambda$
fixing $\dt{-1}$ sends $\alpha_k$ to a multiple of itself. It is
easily seen that if $\lambda$ fixes $d_e$, then it also fixes
$d_{e,a}$.  As a consequence, we can not eliminate the coefficient $a$
by a linear automorphism.  Since the linear part of any automorphism
preserving $d_e$ not only preserves $d_e$ but actually preserves
$d_{c,a}$, it is easy to apply our conjecture to see that the
$d_{c,a}$ are all nonequivalent extensions of $d_c$.  Thus we have
determined all nonequivalent extensions of $\dt{-1}$.
\subsection{The Codifferential $\dstar=\psd{1,1,0 }3$.}\label{0100}
The coboundaries of basic
cochains for \dstar\ are as follows:
\begin{align*}
&D(\phd{1, q, n-q-1}1) =-\psd{1, 1+q, n-q-1}3\\
&D(\phd{0, p, n-p}2) =\psd{1, 1+p, n-p-1}2(n-p)-\psd{1,p, n-p}3\\
&D(\phd{0, p, n-p}3) =\psd{1, 1+p, n-p-1}3(n-p)\\
&D(\psd{0, p, n-p}1) = \phd{1, 1+p, n-p-1}1(n-p)+\phd{0, 1+p, n-p}3\\
&D(\psd{1, q, n-q-1}2) =0\\
&D(\psd{1, q, n-q-1}3) =0
\end{align*}
The cohomology of $\dstar$ is given by
\begin{align*}
&H^1=\langle \psd{1,0,0}2,\psd{1,0,0}3,\phd{0,1,0}3,\phd{0,0,1}3+\phd{1,0,0}1,\phd{0,1,0}2+\phd{0,0,1}3\rangle\\
&H^n=\langle \psd{1,0,n-1}2,\psd{1,0,n-1}3,\phd{0,n,0}2+\phd{0,n-1,1}3,\phd{0,0,n}3+\phd{1,0,n-1}1n\rangle,\quad n>1
\end{align*}

Because the odd part of the cohomology of \dstar\ does not vanish for
$n>2$,  there are nontrivial extensions of this codifferential.
First let us consider when a linear automorphism $\lambda$, given in standard form by \refeq{linauto},
preserves $\dstar$. It is easy to see that $\lambda^*(\dstar)=\dstar$ precisely when $t=0$ and $u=qr$.
Because the cohomology of $d_*$ can be represented by cocycles of type $\delta=\psd{1,0,l}3$
and $\gamma=\psd{1,0,k}2$,  every extension of $\dstar$ is equivalent to one where all added terms are of
one of these two types.
With a little work, one can show that
\begin{align}
\lambda^*(\delta)=&\sum_{x=0}^l\psd{1,x,l-x}3s^xu^{l-x-1}q=\sum_{x=0}^l\psd{1,x,l-x}3\su^x u^{l-1}q
\label{delta}\\
\lambda^*(\gamma)=&
\sum_{x=0}^k\psd{1,x,k-x}2a\su^x\tfrac{qu^k}r-\psd{1,x,k-x}3a\su^{x+1}\tfrac{qu^k}r
\label{gamma}.
\end{align}
We will show that any extension of $\dstar$ is equivalent to an extension of a codifferential of the form
$$d_{k,l}=\psd{1,1,0}3+\psd{1,0,k}2a+\psd{1,0,l}3b.$$
For the moment, we make no assumptions about which of $k$ and $l$ is the larger.
It is usually true that if $a$ and $b$ are nonzero, they can be taken to be 1.
There is one exception to this statement, and that is the case when $k+1=2l$. We will
discuss this special case in more detail later. For the moment, it will be
more convenient for us to leave the coefficients $a$ and $b$ undetermined.

In order to determine the leading terms of coboundaries for the coboundary operator $D_{k,l}$ given by $d_{k,l}$
consider the extended coboundary formulae:
\begin{align*}
D_{k,l}(\phd{1,q,n-q-1}1)&=-\psd{1,1+q,n-q-1}3-\psd{1,q,k+n-q-1}2a-\psd{1,q,l+n-q-1}3b\\
D_{k,l}(\phd{0,p,n-p}2)&=\psd{1,p+1,n-p-1}2(n-p)-\psd{1,p,n-p}3\\&+\psd{1,p-1,k+n-p}2ap+\psd{1,p,l+n-p-1}2b(n-p)\\
D_{k,l}(\phd{0,p,n-p}3)&=\psd{1,p+1,n-p-1}3(n-p)-\psd{1,p,k+n-p-1}2ak
\\&+\psd{1,p-1,k+n-p}3ap+\psd{1,p,l+n-p-1}3b(n-p-l)
\end{align*}

From the coboundary formulas above, we can establish the following \emph{recursion formulas}.
\begin{equation}
D_{k,l}(\phd{0,0,v+1}2)=\psd{1,1,v}2(v+1)-\psd{1,0,v+1}3+\psd{1,0,v+l}2b(v+1)
\end{equation}
\begin{multline}
D_{k,l}(\phd{0,p,v+1}3-\phd{1,p,v}1k)=
\psd{1,p+1,v}3(v+k+1)
+\psd{1,p-1,v+k+1}3ap\\
+\psd{1,p,v+l}3b(v+k+1-l)
\end{multline}
\begin{multline}
D_{k,l}(\phd{0,p+1,v+1}2-\phd{1,p,v+1}1+\phd{0,p,v+l+1}2b)=\\
\psd{1,p+2,v}2(v+1)
+\psd{1,p+1,v+l}2b(l+2+2v)
+\psd{1,p,v+k+1}2a(p+2)\\
+\psd{1,p,v+2l}2b^2(l+1+v)
+\psd{1,p-1,v+k+1+l}2abp
\end{multline}
If $\eta-\xi$ is a $D_{k,l}$ coboundary, then let us denote this by $\eta\sim\xi$.
The recursion formulas above allow us to conclude the following \emph{reduction formulas}.
\begin{align}\label{reduction}
\psd{1,p+1,v}3&\sim-\psd{1,p,v+l}3b\tfrac{v+k+1-l}{v+k+1}-\psd{1,p-1,v+k+1}3a\tfrac{p}{v+k+1}\\
\psd{1,p+2,v}2&\sim -\psd{1,p+1,v+l}2b\tfrac{l+2(v+1)}{v+1}-\psd{1,p,v+k+1}2a\tfrac{p+2}{v+1}
\\&-\psd{1,p,v+2l}2b^2\tfrac{l+v+1}{v+1}-\psd{1,p-1,v+k+1+l}2ab\tfrac{p}{v+1}\nonumber\\
\psd{1,1,v}2&\sim\psd{1,0,v+1}3\tfrac1{v+1}-\psd{1,0,v+l}2b
\end{align}
The first reduction formula does not hold when $a=0$; instead, we have the simpler reduction formula
\begin{equation}
\psd{1,p+1,v}3\sim-\psd{1,p,v+l}3b.
\end{equation}
These formulas show us to reduce any cochain of the form $\psd{0,p,n-p-1}2$ or $\psd{0,p,n-p-1}3$
to a cochain where the middle index of each term is zero, modulo a coboundary.

If we are considering an extension of $d_{k,l}$, we know that it can be reduced to one of the form
\begin{equation*}
d=\psd{1,1,0}3+\psd{1,0,k}2a+\psd{1,0,l}3b +\psd{1,0,m}2a_m+\psd{1,0,m}3b_m,
\end{equation*}
where $a_m=0$ if $m\le k$ and $b_m=0$ if $m\le l$.  Thus we really are interested in when we can get
rid of terms of the form $\psd{1,0,m}2$ and $\psd{1,0,m}3$ in the expression above,
which can be done only when they appear as leading terms in $D_{k,l}$ coboundaries.   The reduction formulas allow us to
add coboundary terms to a coboundary to reduce terms of the form $\psd{1,q,n-q-1}2$ and
$\psd{1,q,n-q-1}3$ to terms of the form $\psd{1,0,m}2$ and $\psd{1,0,m}3$, for certain values
of $m$.  This will help us to determine which terms of this form are leading terms in coboundaries.

We will be studying the coboundaries of the even cochains
\begin{align*}
&\ph_n=\phd{1,0,n-1}1n+\phd{0,0,n}3\\
&\ph'_n=\phd{0,n,0}2+\phd{0,n-1,1}3,
\end{align*}
which are a basis of the even part of $H^n(\dstar)$.  What we expect is that these $D_*$-cocycles will give rise to
new $D_{k,l}$-coboundaries, which will allow us to eliminate certain terms in the expressions above
which were not leading terms for $D_*$-coboundaries, but which are leading terms for $D_{k,l}$-coboundaries.
We have
\begin{align}
D_{k,l}(\ph_n)=&-\psd{1,0,k+n-1}2a(k+n)-\psd{1,0,l+n-1}3bl\label{Dphin}\\
D_{k,l}(\ph'_n)=&-\psd{1,n-1,k}2a(k-n)+\psd{1,n-2,k+1}3a(n-1)\label{Dphinprime}\\&+\psd{1,n-1,l}3b(1-l)\nonumber.
\end{align}

Two extensions of $\dstar$ can only be equivalent if $\lambda^*$ applied to the secondary term of the first
extension differs from the secondary term in the the second extension by a coboundary for some linear automorphism
$\lambda$.  Thus we first need to consider the action of the linear automorphism group on the $D_*$-cohomology.
The secondary term can be taken to be of the form $\psd{1,0,k}2a+\psd{1,0,k}3b$, because every $D_*$-cohomology
class of degree $k$ can be represented by a $D_*$-cocycle of this form. If $a$ or $b$ vanish, then applying a
diagonal automorphism one sees easily that the other coefficient can be taken to be 1,
and similarly, if both do not vanish, they can both be taken to be 1. Now
\begin{equation*}
\lambda^*(\psd{1,0,k}2+\psd{1,0,k}3)=
\sum_{x=0}^k\psd{1,x,k-x}2\su^x\tfrac{qu^k}r+\psd{1,0,k}3\su^x(u^{k-1}r-\su\tfrac{qu^k}r).
\end{equation*}
If you choose $\lambda$ so that $u=qr$, ${qu^k}r=1$ and $\su=u^{k-1}r$, then the above cocycle is cohomologous
to $\psd{1,0,k}2$.  As a consequence,  we can assume that the secondary term is either of the
form $\psd{1,0,k}2$ or $\psd{1,0,l}3$ for some $k$ or $l$.
\subsubsection{Extensions of $\dstar$ with secondary term $\psd{1,0,k}2$}\label{Remark 1}
Let
$$
d_e=\psd{1,1,0}3+\psd{1,0,k}2.
$$
Recall that for any generalized automorphism of $d_e$, its linear part $\lambda$ preserves $\dstar$
and $\lambda^*(\psd{1,0,k}2)$, given by \refeq{gamma}, is $D_*$-cohomologous to $\psd{1,0,k}2$. It follows that
$u=qr$, $qu^k=r$ and $s=0$.
But this implies that $\lambda^*(d_e)=d_e$. This fact greatly simplifies the study of extensions
of $d_e$.

The $D_e$-coboundaries of the $D_*$-cohomology classes are given by
\begin{align*}
D_e(\ph_n)=&-\psd{1,0,k+n-1}2(k+n)\\
D_e(\ph'_n)=&-\psd{1,n-1,k}2(k-n)+\psd{1,n-2,k+1}3(n-1).
\end{align*}

Thus $\psd{1,0,m}2$ is always a $D_e$-coboundary if $m\ge k$.
The case with $\ph'_n$ is more complicated.
Note that in the first two
reduction formulas, since $b=0$, the middle upper index on the right hand side always drops by
2, and the right hand side of the first reduction formula
vanishes for $\psd{1,1,v}3$. Thus we conclude that any term of the form
$\psd{1,2m+1,v}3$ is a $D_e$-coboundary.  Moreover, any term of the form $\psd{1,2m,v}2$ reduces to a multiple
of $\psd{1,0,v+m(k+1)}2$, so is also a $D_e$-coboundary.
When $n=2m+1$ is odd,
$$D_e(\ph'_n)=\psd{1,2m,k}2(2m+1-k)+\psd{1,2m-1,k+1}3(2m),$$
which is a $D_e$-coboundary, so $\ph'_{2m+1}$
extends to a $D_e$-cocycle.

Notice that for $d_e$, the reduction formulas simplify to
\begin{align*}
\psd{1,p+1,v}3\sim&-\psd{1,p-1,v+k+1}3\tfrac p{v+k+1}\\
\psd{1,p+2,v}2\sim&-\psd{1,p,v+k+1}2\tfrac{p+2}{v+1}\\
\psd{1,1,v}2\sim&\psd{1,0,v+1}3\tfrac 1{v+1}
\end{align*}

When $n=2m$, applying the reduction formulas to the coboundary of $\ph'_n$
yields a coboundary of the form $\phd{1,0,m(k+1)}3$, multiplied by a nonzero constant.
Since
$$D_e(\ph'_{2m})=\psd{1,2(m-1)+1,k}2(2m-k)+\psd{1,2(m-1),k+1}3(2m-1),$$
we note that after applying the recursion formulas, both terms reduce to multiples of
$\psd{1,0,m(k+1)}3$. Adding the coefficient arising in reducing the first term
to the one arising in reducing the second term gives
an overall nonzero coefficient
\begin{equation}\label{prodform}
\frac{(-1)^{m-1}\prod_{i=1}^{m}(2i+1)}{(k+1)^m m!}
\end{equation}
for the leading cochain $\psd{1,0,m(k+1)}3$ in the reduced form of $D_e(\ph'_{2m})$.

Thus,  $\ph'_{2m}$ gives rise to a nontrivial $D_e$-coboundary, and it does not extend
to a $D_e$-cocycle. The even part of the cohomology of $d_e$ has a basis given by the
extensions of  $\ph'_{2m+1}$ to cocycles,
while the odd part of the cohomology has a basis given by $\psd{1,0,l}3$ for those $l$ which are not
multiples of $k+1$.
As a consequence, in an extension of $d_e$ of the form $d_{k,l}+\ho$, with $k<l$, if $l$ is a multiple of $k+1$,
we can add a $D_e$-coboundary to eliminate this term.  Thus, in classifying extensions of $d_e$, we
only need consider those for which $l$ is not a multiple of $k+1$. In particular,
we do not need to consider the case when $l=k+1$.
In section \ref{Remark 3} we will consider extensions of $d_e$ of the form
$$
d_{k,l}=\psd{1,1,0}3+\psd{1,0,k}2+\psd{1,0,l}3,\qquad l\ne m(k+1).
$$
The reduction formulas of the previous section will enable us to determine the leading terms of $D_{k,l}$-coboundaries.
We will also have occasion to consider extensions of $d_{k,l}$ of the form
\begin{equation}
d_{k,l,m}=\psd{1,1,0}3+\psd{1,0,k}2+\psd{1,0,l}3+\psd{1,0,m}3b,
\end{equation}
and will want recursion and reduction formulas for this codifferential as well. The modified recursion formulae are
\begin{multline}\label{exrec1}
D_{k,l,m}(\phd{0,0,v+1}2)=\\\psd{1,1,v}2(v+1)-\psd{1,0,v+1}3+\psd{1,0,v+l}2(v+1)+\psd{1,0,v+m}2b(v+1)
\end{multline}
\begin{multline}\label{exrec2}
D_{k,l,m}(\phd{0,p,v+1}3-\phd{1,p,v}1k)=
\psd{1,p+1,v}3(v+k+1)
+\psd{1,p-1,v+k+1}3p\\
+\psd{1,p,v+l}3(v+k+1-l)
+\psd{1,p,v+m}3b(v+k+1-m)
\end{multline}
\begin{multline}\label{exrec3}
D_{k,l,m}(\phd{0,p+1,v+1}2-\phd{1,p,v+1}1+\phd{0,p,v+l+1}2+\phd{0,p,v+m+1}2b)=\\
\psd{1,p+2,v}2(v+1)
+\psd{1,p+1,v+l}2(l+2+2v)+\psd{1,p+1,v+m}2b(m+2+2v)\\
+\psd{1,p,v+k+1}2(p+2)+
+\psd{1,p,v+2l}2(l+1+v)+\psd{1,p,v+2m,}2b^2(m+1+v)\\
+\psd{1,p-1,v+k+1+l}2p+\psd{1,p-1,v+k+1+m}2bp
\end{multline}

\subsubsection{Extensions of $\dstar$ with secondary term $\psd{1,0,l}3$}\label{Remark 2}
Let
$$d_e=\psd{1,1,0}3+\psd{1,0,l}3.$$
The linear part $\lambda$ of any generalized automorphism of $d_e$ preserves $\dstar$, and so its action on
$\delta=\psd{1,0,l}3$ is given by \refeq{delta}. Since we must have $\lambda^*{\delta}\sim \delta$, it follows
that $u^{l-1}q=1$, and this is the only condition necessary.  Of course, when $\lambda^*(\delta)\ne\delta$, we
must follow $\lambda^*$ with some formal automorphism $g^*$ such that $g^*(\lambda(\delta))=\delta +\ho$, in order
to obtain a generalized automorphism $\lambda g$ of $d_e$. In fact, if we choose
$$\alpha=-\sum_{x=0}^l\phd{0,x-1,l-x+1}3\tfrac 1{l-x+1}\bigl(\tfrac su \bigr)^x,$$
then $\exp(-\ad\alpha)(\lambda(\delta))=\delta+\ho$.

Because $D_e(\ph_n)=-\psd{1,0,l+n-1}3l$, we see that $\psd{1,0,m}3$ occurs as the leading term of a $D_e$-coboundary
for any $m\ge l$, so
in any extension of $d_e$, we can assume that no such terms occur. Thus, up to equivalence, any nontrivial extension
of $d_e$ is of the standard form
$$
d_{k,l,e}=\psd{1,1,0}3+\psd{1,0,l}3+\psd{1,0,k}2a +\ho,
$$
where all higher order terms are of the form $\psd{1,0,m}2a_m$, and $a\ne0$. Since the linear part $\lambda$ of a generalized
automorphism $\lambda g$ of $d_{k,l,e}$ is the linear part of a generalized automorphism of $d_e$, we know that
$u^{l-1}q=1$, and $u=qr$. Note that if $2l\ne k+1$, then from \refeq{gamma}, we can find a diagonal automorphism $\lambda$
such that the coefficient of $\psd{1,0,k}2$ in $\lambda^*(d_{k,l})$ is 1. Moreover, in that case, then any generalized
automorphism of $d_{k,l}$ will satisfy $u^{k+1-2l}=1$, which limits the diagonal part of $\lambda$ to just a few
possibilities.

Let us consider a formal automorphism $g=\exp(-\alpha_2)\cdots$ such that $\lambda g$ is a generalized
automorphism of $d_{k,l}$. Then
$Q=g^*\lambda^*(d_{k,l})$ is of the form $d_{k,l}+\ho$. We will show that for a certain value of $m$ the coefficient
of $\psd{1,0,m}2$ appearing in $Q$ is exactly $a_m$ plus a nonzero multiple of $s$. Thus, for an appropriate value
of $s$, the coefficient becomes zero. As a consequence, if we only consider extensions of $d_{k,l}$ such that the
coefficient of $\psd{1,0,m}2$ is zero,  we do not lose any of the equivalence classes.  Moreover, if an automorphism
of $d_{k,l,e}$ preserves this property, then the coefficient $s$ in its linear part
must vanish. This observation will allow us to restrict
our consideration to automorphisms whose linear part is diagonal.
In most cases, there are
only a few diagonal matrices satisfying our requirements, and they always
preserve $d_{k,l}$. Because of this, we will be able
to use our main conjecture to classify the extensions of $d_{k,l}$.

If $n<l$, then $\alpha_n$ must be a $D_*$-cocycle; otherwise terms of degree lower than $l$ would appear
in $Q$. Moreover, no term of
type $\phd{0,p,n-p}2$ can appear in $\alpha_n$ for any $n<k$, because otherwise $Q$ would contain
a term of type $\psd{1,p+1,n-p-1}2$. If $n<l$, the only combination of
$\phd{0,0,n}3$
and $\phd{1,0,n-1}1$ which can appear in $\alpha_n$ is a multiple of $\ph_n$, since this is the
only $D_*$-cocycle combination of these terms.
However, since $d_e(\ph_n)=-\psd{1,0,l+n-1}3l$, and there
are no terms like this in $Q$ if $1<n<k-l+1$, the coefficient of $\ph_n$ must vanish unless
$n\ge k-l+1$.

If $l\le n< k-l+1$ and
terms of type $\phd{0,0,n}3$ and $\psd{1,0,n-1}3$ occur in $\alpha_n$, then since
$d_e(\phd{0,0,n}3)$ and $d_e(\psd{1,0,n-1}3)$ each contain a term of type $\psd{1,0,l+n-1}3$,
then they appear as a multiple of the cochain
$\eta_n^0$, where $$\eta_n^p=\phd{1,p,n-p-1}1\pfrac{n-p-l}l +\phd{0,p,n-p}3\pfrac 1l.$$
Since $\eta_n^0$ and $\ph_n$ span the two dimensional
subspace spanned by $\phd{0,0,n}3$ and $\phd{1,0,n-1}1$, when $n\ge k-l+1$ we can
discuss their contributions to $\alpha_n$ separately.

Define
$$
\zeta_n^p=\phd{1,p,n-p-1}1(n-p)+\phd{0,p,n-p}3 \quad\text{if $0\le p\le n-1$}.
$$
Then
$$\alpha_l=-\sum_{p=0}^{l-1}\eta_l^p \su^{p+1}+\zeta_n^p c_l^p,$$
where $c_n^p$ are arbitrary constants, and $c_l^0=0$.
In computing $g^*\lambda^*(d_{k,l,e})$, we will be computing terms of the form
$(\ad\alpha_j)^kD^\lambda_{k,l,e}(\alpha_i)$ where $D^\lambda_{k,l}$ is the coboundary
operator determined by $\lambda^*(d_{k,l,e})$.
Now
\begin{align*}
D^\lambda_{k,l}(\alpha_l)=&
\sum_{p=0}^{l-1}
\Bigl[-\psd{1,p+1,l-p-1}3\su^{p+1}\\
&+\sum_{x=0}^l\psd{1,p+x,2l-x-p-1}3\su^x(c_l^p(x-l)-\su^{p+1}\tfrac xl)\\
&+\sum_{v=k}^\infty\sum_{x=0}^va_v\bigl[
\\&+\psd{1,p+x,v+l-x-p-1}2\su^x({ c_l^p(x+p-v-l)-\su^{p+1}\tfrac{x+p-v}l})
\\&+\psd{1,x+p-1,v+l-p-x}3\su^{x}(pc_l^p-\su^{p+1}\tfrac{p}l)
\\&+\psd{1,x+p,v+l-p-x-1}3\su^{x+1}(c_l^p(v-x)+\su^{p+1}\tfrac{x+l-v}l)
\bigr]\Bigr]
\end{align*}
The first set of terms have degree $l+1$ and
are necessary to eliminate the terms of type $\psd{1,p,l-p}3$ from $\lambda^*(\delta)$,
the second set have
degree $2l$, while all the rest have degree $v+l$. Notice that in the coboundary above,
the coefficient of the $\psd{1,0,2l}3$ term
is zero, the coefficient of the $\psd{1,0,k+l-1}2$ term is $\su\tfrac kl$, while the coefficient of the
$\psd{1,0,k+l-1}3$ depends on the choice of the coefficients $c_l^p$.

Next,  remember that in applying an exponential of a coderivation,  we also obtain terms
of the form $(\ad\alpha_l)^i(D^\lambda_{k,l}(\alpha_l)$.  Let us examine the first such term.
$[\alpha_l,\D_{k,l}(\alpha_l)]$ will contain some terms of degree $2l$, $3l-1$, and $v+2l-1$. The terms of
degree $2l$ and $3l-1$ are of the type $\psd{1,p,x-p}3$ where $p>0$. No term of type $\psd{1,0,x}3$ of
or any term of type $\psd{1,p,x-p}2$ of degree less than $k+2l-1$ will
arise in $(\ad\alpha_l)^i(D^\lambda_{k,l}(\alpha_l))$ for $i\ge 1$.

For $n<l$, $\alpha_n$ consists only of multiples of $\eta_n^p$, with $p>0$, unless $n\ge k+l-1$,
in which case, note that $\ph_n=\eta_n^0$, so we will discuss these terms when we discuss the
contribution of the $\ph_n$ terms.

We have addressed how to construct $\alpha_n$ in order to remove all unwanted terms of type $\psd{1,p,n-p}3$
for $p>0$.  The first place we encounter terms that must be eliminated which are not of this form is in
degree $k+1$ where they come from $\lambda^*(\gamma)$. The terms in $\lambda^*(\gamma)$ which concern
us are those of type $\psd{1,0,k}3$ and $\psd{1,1,k-1}2$, the latter because it gives rise to a term of the
former type in the process of elimination. We have
\begin{multline*}
D^\lambda_{k,l,e}(-\phd{0,0,k}2\su\tfrac ak)=
-\psd{1,1,k-1}2\su a+\psd{1,0,k}3\su \tfrac ak\\
+\sum_{x=0}^k(-\psd{1,x,k+l-x-1}2\su^{x+1}a+\psd{1,x-1,k+l-x}3\su^{x+1}\tfrac{ax}k)+\ho
\end{multline*}
The only terms that are important to us in the expression above are the first two, and
$-\psd{1,0,k+l-1}2\su a$.
To get rid of the the $\psd{1,0,k}3$ term, we use $\ph_{k-l+1}$. We have to remember that in addition to
the coefficient $-a\su$ arising from $\lambda^*(\gamma)$, we also must remove the coefficient $\pfrac ak\su$
which arises above. The total coefficient is $\su\tfrac{a(1-k)}k$. We really only need to look at
\begin{equation*}
D_{k,l}\bigl(\ph_{k-l+1}\su\tfrac{a(1-k)}{kl}\bigr)=
-\psd{1,0,k}3\su\tfrac{a(1-k)}{k}-\psd{1,0,2k-l}2\su\tfrac{a^2(1-k)(2k-l+1)}{kl}.
\end{equation*}
When $k+1<2l$, it follows that $k<2k-l<k+l-1$, and the
only occurrence of $\psd{1,0,2k-l}2$ in the calculation of the addition of coboundaries and higher order
terms coming from the exponential is from the term above. Thus, by an appropriate choice of $s$, we can arrange that a
standard form of an equivalent codifferential will have a zero value for the coefficient.

When $k+1>2l$, it follows that $k+l-1<2k-l$, and in addition to the contribution to the $\psd{1,0,k+l-1}2$ term from the
above coboundary, we also have a contribution from the $\eta_l^0$ coboundary. The overall coefficient added is
$\su\tfrac{a(k-l)}l$.  Thus we can reduce the $\psd{1,0,k+l-1}2$ to zero by an appropriate choice of $s$.

Finally, when $k+1=2l$, the terms of type $\psd{1,0,2k-l}2$ and $\psd{1,0,k+l-1}2$
are both of type $\psd{1,0,2l-1}2$, with total coefficient
\begin{math}
\su\tfrac{a(-1+l)(6al-2a+2l-1)}{(2l-1)l},
\end{math}
which is zero precisely when $a=-\tfrac{2l-1}{2(3l-1)}$.
Otherwise, the same idea works, and we can again assume that the $\lambda$ is
diagonal and that our standard form of the codifferential has a zero coefficient for the $k+l-1$ spot. We
will discuss this case in more detail in the subsection devoted to the special case $k+1=2l$.

Now let us make a few remarks about $d_e$ and its cohomology.
Equation (\ref{Dphin}) reduces to
$$\psd{1,0,m}3=D_e(-\ph_{m-l+1}\tfrac1l),\quad\text{if } m\ge l.$$
Since the $D_e$-coboundary of $\ph'_n$ can be reduced to an element of this same form, it
follows that $\ph'_n$ can be extended to a $D_e$-cocycle. The exception is that in the case $n=1$,
the term we add to $\ph'_1$ is of the same degree, so, instead of extending to a cocycle, we
see that $$\ph_1(1-l)+\ph'_1l$$ is a $D_e$-cocycle. Moreover, $D_e(\phd{0,1,0}3)=-\psd{1,1,l-1}3$,
which is a $D_e$-coboundary, so $\phd{0,1,0}3$ extends to a $D_e$-cocycle.
Thus, we obtain the following table for the cohomology of $d_e$.
\begin{align*}
&H^1=\langle \psd{1,0,0}2,\psd{1,0,0}3,\phd{0,1,0}3+\ho,\ph_1(1-l)+\ph'_1l\rangle\\
&H^n=\langle \psd{1,0,n-1}2,\ph'_n+\ho\rangle,\quad n>1
\end{align*}

Occasionally we will need an expanded version of the reduction formulas, related to an extension
of $d_{k,l}$ of the form
\begin{equation}
d_{k,l,m}=
\psd{1,1,0}3+\psd{1,0,l}3+\psd{1,0,k}2a+\psd{1,0,m}2b.
\end{equation}
The main difference in the nature of the recursion formulas for $d_{k,l,m}$ is that the second formula
has to be applied to an infinite sum $\theta_{p,v}$ of cochains defined as follows:
\begin{equation}
\theta_{p,v}=
\psd{0,p,v+1}3-\psd{1,p,v}1k+\sum_{i=1}^\infty\psd{1,p,v+i(m-k)}1(-1)^i(b/a)^i(m-k).
\end{equation}
We calculate
\begin{multline*}
D_{k,l,m}(\theta_{p,v})=
\psd{1,p+1,v}3(v+1+k)+\sum_{i=1}^\infty\psd{1,p+1,v+i(m-k)}3(-1)^{i+1}(b/a)^i(m-k)\\
+\psd{1,p,v+l}3(k+v+1-l)+\sum_{i=1}^\infty\psd{1,p,v+l+i(m-k)}3(-1)^{i+1}(b/a)^i(m-k)\\
+\psd{1,p-1,v+k+1}3pa+\psd{1,p-1,v+m+1}3pb.
\end{multline*}
In order to obtain a recursion formula from the formula above, we need to reduce the
number of occurrences of the $p+1$ index to 1. Accordingly, we define
\begin{equation*}
\Theta_{p,v}=\theta_{p,v}+\sum_{i=1}^\infty\theta_{p,v+i(m-k)}(-1)^i\pfrac{(b/a)^i(m-k)}{v+m+1}
\end{equation*}
The modified recursion formulas are as follows.
\begin{equation}\label{recur1}
D_{k,l,m}(\phd{0,0,v+1}2)=\psd{1,1,v}2(v+1)-\psd{1,0,v+1}3+\psd{1,0,v+l}2(v+1).
\end{equation}
\begin{multline}\label{recur2}
D_{k,l,m}(\Theta_{p,v})=
\psd{1,p+1,v}3(v+1+k)
\\+\psd{1,p,v+l}3(k+v+1-l)+\sum_{i=1}^\infty\psd{1,p,v+l+i(m-k)}3\pfrac{(-1)^{i+1}(b/a)^i(m-k)l}{v+m+1}\\
+\psd{1,p-1,v+k+1}3pa+\psd{1,p-1,v+m+1}3pb\pfrac{v+k+1}{v+m+1}.
\end{multline}
\begin{multline}\label{recur3}
D_{k,l,m}(\psd{0,p+1,v+1}2-\psd{1,p,v+1}1+\psd{0,p,v+l+1}2)=\\
\psd{1,p+2,v}2(v+1)+\psd{1,p+1,v+l}2(l+2+2v)\\
+\psd{1,p,v+2l}2(l+v+1)+\psd{1,p,v+k+1}2a(p+2)+\psd{1,p,v+m+1}2b(p+2)\\
+\psd{1,p-1,v+k+l+1}2ap+\psd{1,p-1,v+m+l+1}2bp.
\end{multline}

\subsubsection{General Extensions of $\dstar$, with $2l\ne k+1$}\label{Remark 3}
We consider the codifferential
\begin{equation*}
d_{k,l}=\psd{1,1,0}3+\psd{1,0,k}2+\psd{1,0,l}3
\end{equation*}
When $k<l$, this codifferential arises as an extension of the codifferential considered in (\ref{Remark 1}),
so we can assume that $l$ is not a multiple of $k+1$. When $l< k$, this codifferential
arises as an extension of the codifferential considered in (\ref{Remark 2}), so will will only
consider extensions of $d_{k,l}$ in which the coefficient of $\psd{1,0,k+l-1}2$ vanishes.

No matter whether this is an extension of the $d_e$ in (\ref{Remark 1}) or the one in (\ref{Remark 2}),
the condition that $\lambda^*(d_{k,l})-d_{k,l}$ be the leading term of a $d_e$-coboundary is satisfied
precisely when $u=qr$, $r=qu^k$ and $u^{l-1}q=1$,  so we can solve $r=u^l$, $q=u^{1-l}$, $u^{2l-(k+1)}=1$.
Moreover,
it is easily checked that $\lambda^*(d_{k,l})=d_{k,l}$. There are only at most $|2l-(k+1)|$
solutions for $\lambda$, and $\lambda$ applied to a cochain of fixed degree simply multiplies it by
a power of $u$.

The group $U_{k,l}$ of $|2l-(k+1)|$ roots of unity acts in an obvious way on the
set of extensions of $d_{k,l}$.
Thus, in studying extensions of
$d_{k,l}$, we can determine equivalence by studying coboundaries.  Our main goal in the following will
be to add terms to $d_{k,l}$ until all of the even cocycles have been killed off. Once we have arrived
at an extension $\dfin$ of $d_{k,l}$ for which no even cohomology classes remain, the set of equivalence classes
of $\dfin$ is given by the set of extensions of $\dfin$ such that no term of degree larger than the maximal
degree in $\dfin$ is a coboundary with respect to $\dfin$, modulo the action of the group $U_{k,l}$.

The formula for the coboundary of $\ph_n$ allows us to convert cocycles of the form $\psd{1,0,m}2$ to
the form $\psd{1,0,m+l-k}3$ and vice versa, up to a coboundary term.
Explicitly, we have the conversion formulas
\begin{align*}
\psd{1,0,v}2&\sim-\psd{1,0,v+l-k}3\tfrac l{v+1}\\
\psd{1,0,v}3&\sim-\psd{1,0,v+k-l}2\tfrac {v+k+1-l}{l}
\end{align*}
Studying the reduction formulas, observe that
whenever a middle upper index reduces by 2 in a reduction formula, the right upper index increases
by either $k+1$ or $2l$.
Whenever the middle upper index reduces by 1 in a reduction formula, the right upper index increases by
$l$, with the exception that for $\psd{1,1,v}2$, the right upper index in the first term
only increases by 1, but in this case, we also convert the lower index from a 2 to a 3,
so if you convert the term with a lower index 3 back into a term with lower index 2, the net
effect is that the right upper index increased by $k+1-l$. At the last step,
we also can convert terms with lower index 2 to terms with lower index 3, and then the upper right hand
coefficient increases by $l-k$.

Looking at the first term in $D_{k,l}(\ph'_n)$, and noting that its upper right index is $k$, we see that
in reductions to terms of middle upper index 0 and lower index 3,
at the last step, the $k$ is replaced by an $l$. Putting all these facts
together we obtain that the reduced form of $D_{k,l}$ contains only terms of the type
$\psd{1,0,i(k+1)+jl}3$, where $2i+j=n$, and every such term arises in the
reduction process, possibly with a net zero coefficient. Moreover, if $2i+j=2i'+j'$, then
$i(k+1)+jl\ne i'(k+1)+j'l$, since $2l\ne k+1$.

If $k+1<2l$, then if $n=2m$, the smallest upper right index  is $m(k+1)$, while if $n=2m+1$, then the smallest
such index is $m(k+1)+l$.
If $k+1>2l$, then the smallest
upper right index is always just $nl$. Notice that when $k+1=2l$ all coefficients have the same index. Thus
the three cases are best treated separately. Note that in all three cases, the smallest upper right index appearing
in the expression is different for different values of $n$, so $\ph'(n)$ and $\ph'(n')$ do not reduce to terms
with the same order, if the coefficient of the smallest degree term is nonzero.

\centerline{\emph{\bf Case 1: $k+1<2l$}}
The analysis depends somewhat on whether $n$ is even or odd, so we treat these cases separately.
If $n$ is even, say $n=2m$, then the coefficient of the $\psd{1,0,m(k+1)}3$ term is the
the same coefficient given in \refeq{prodform}.
In particular, it is nonzero. When $n$ is odd, the situation is more complicated.

If $n=2m+1$, the coefficient of the $\psd{1,0,m(k+1)+l}3$ term is
$$
\frac{(-1)^m(-2l+(2m+1)(k+1))\prod_{i=1}^m(2l+(2i+1)(k+1))}{(k+1)^{m+1}\prod_{i=1}^m (l+i(k+1))},
$$
which vanishes only when $2l=n(k+1)$. Thus if $k+1<2l$, we can only have a zero coefficient for at most one
value of $n$, and for most values of $k$ and $l$ this doesn't occur.
Since there are some differences in these cases, we treat them separately.

\centerline{\emph{ Subcase 1: $n=2m+1$ and  $2l=n(k+1)$}}
Note that this case only occurs when $k<l$, since $k+1<2l$.
Note that if $x\ge 0$, then
\begin{equation}i(k+1)+(2x+1)l=(i+nx)(k+1)+l.\end{equation}
If $2i+j=n$ and $j>1$, then $j=2x+1$ for some $x>0$. But then the term of type
$\psd{1,0,i(k+1)+jl}3$ in $D_{k,l}(\ph'_n)$ is the same type as the leading term
of $D_{k,l}(\ph'_{2(i+nx)+1})$. Since the coefficient of the term $\psd{1,0,m(k+1)+l}3$ vanishes,
this shows that $D_{k,l}\ph'_n$ is a sum of coboundaries of different $\ph's$,

Moreover, every upper right index in any nonleading term in
the coboundary of any of the $\ph'$ terms occurs as the upper
right index of a leading term of some other such coboundary.  As a consequence,
every term in the coboundary of any $\ph'$, with the exception of its leading order term, occurs
as a leading order term of the coboundary of some $\ph'_{n'}$, with $n\ne n'$
Thus, by subtracting appropriate coboundaries
of higher order $\ph'$s, one arrives at the conclusion that $\phd{1,0,m'(k+1)}3$ and $\phd{1,0,m'(k+1)+l}3$
are actual coboundaries, not just leading order terms of coboundaries,  with the exception of $\phd{1,0,m(k+1)+l}3$.

Moreover,  $\ph'_{2m+1}$ extends to a $D_{k,l}$-cocycle.  Thus exactly one higher order even cohomology class
remains
and many odd cohomology classes remain, those whose upper right index is not a multiple of $(k+1)$ or
of the form $m'(k+1)+l$
for some $m'\ne m$.
In order to complete the classification of our extensions,  we need to go one step further,  because we have
not yet killed off all the even cohomology.  Consider an extension
\begin{equation}
d_{k,l,x}=\psd{1,1,0}3+\psd{1,0,k}2+\psd{1,0,l}3+\psd{1,0,x}3b
\end{equation}
of $d_{k,l}$, where we assume that $x$ is not a multiple of $k+1$, nor is
$x$ of the form  $x=y(k+1)+l$ for any $y\ne m$. (Otherwise, the extension is equivalent to $d_{k,l}$.)
Using the extended recursion formulas (\ref{exrec1},\ref{exrec2},\ref{exrec3}) from section \ref{Remark 1},
we can compute that
the coefficient of the term of type $\psd{1,0,m(k+1)+x}3$ term in the reduced form of $D_{k,l,m}(\ph'_n)$
is
\begin{equation}
\frac{(-1)^m(-2x+(2m+1)(k+1))\prod_{i=1}^m(2x+(2i+1)(k+1))}{(k+1)^{m+1}\prod_{i=1}^m (x+i(k+1))},
\end{equation}
which does not vanish.  Moreover, our condition on $x$ guarantees that the index $m(k+1)+x$ is not of the
form $v(k+1)+l$.  Since the upper right index in the leading term in the reduced form of
$D_{k,l,m}(\ph'_{n'})$ is of the form $v(k+1)+l$ when $n'\ne n$, we know that $m(k+1)+x$ is not the upper
right index of a $D_{k,l}$-coboundary. Of course, there are terms in the reduced form of $D_{k,l,m}(\ph'_n)$
which are of the form $i(k+1)+jl$, where $2i+j=n$, which may have smaller degree than $m(k+1)+x$. However,
as  before, they occur as leading terms of $D_{k,l,m}$-coboundaries of $\ph's$ of higher degree,
which have reduced forms with indices with  $k+1$ and $l$ terms (which can be expressed in terms of
even higher degree $\ph's$), and terms with indices involving $k+1$ and $l$ and $x$ of degree higher than
$m(k+1) +x$.  Thus, we are able to conclude that $\psd{1,0,m(k+1)+x}3$ is the leading term of the reduced form
of the coboundary of $\ph'_n$.  Thus we have finally killed off all the even cohomology.
Any two extensions of $d_{k,l,x}$ of the form
\begin{equation}
d_{k,l,x,e}=\psd{1,1,0}3+\psd{1,0,k}2+\psd{1,0,l}3+\psd{1,0,x}3b+\sum_{y=x+1}^\infty\psd{1,0,y}3b_y,
\end{equation}
where $b_y=0$ if $y$ is a multiple of $k+1$, is of the form $y=z(k+1)+l$ for some $z\ne m$, or is equal
to $m(k+1)+x$, are equivalent precisely when they are equivalent under the action of the group $U_{k,l}$
of $2l-(k+1)$ roots of unity.

\centerline{\emph{ Subcase 2: $n=2m+1$ and $2l\ne n(k+1)$}}

Note that when $k\le l$, the leading terms of the coboundaries associated to $\ph_n$ and $\ph'_n$ are
$\psd{1,0,k+n-1}2$ and $\psd{1,0,m(k+1)}3$ or $\psd{1,0,m(k+1)+l}3$ depending on whether $n$ is even or odd.
Thus, when $k\le l$ we obtain that any extension of $d_{k,l}$ is equivalent to one of the form

\begin{equation}
d_{k,l,x}=\psd{1,1,0}3+\psd{1,0,k}2+\psd{1,0,l}3+\sum_{x=l+1}^\infty\psd{1,0,x}3b_x,
\end{equation}
where $b_x=0$ if $x$ is a multiple of $k+1$, or $l$ plus a multiple of $k+1$.

When $l<k+1<2l$, then the terms we should be adding should be converted to the opposite kind.
Thus the leading terms of coboundaries are $\psd{1,0,l+n-1}3$ and $\psd{1,0,(m+1)(k+1)-l-1}2$ or
$\psd{1,0,(m+1)(k+1)-1}2$, depending on whether $n$ is even or odd.
Thus any extension of $d_{k,l}$ is equivalent to one of the form
$$
d_{k,l,e}=\psd{1,1,0}3+\psd{1,0,l}3+\psd{1,0,k}2+\psd{1,0,x}2a_x,
$$
where $a_x=0$ if $x=(m+1)(k+1)-l-1$ or $x=(m+1)(k+1)-1$, and this condition classifies these extensions
up to the action of the group $U_{k,l}$.

\centerline{\emph{\bf Case 2: $k+1>2l$}}

In this case, the lowest upper right index of $D_{k,l}(\ph'_n)$ is $nl$ and
the coefficient of this term  is simply
\begin{equation}
\frac{(-1)^n((n+1)l-(k+1))}{(n-1)l+k+1},
\end{equation}
which only vanishes if $k+1=nl$. Again, we need to treat this
special case separately.

\centerline{\emph{ Subcase 1: $k+1=(n+1)l$}}

This case is similar to the case when $2l=n(k+1)$ above, in that
$\psd{1,0,ml}3$ is a coboundary, not just
the leading term of a coboundary, when $m\ne n$. Moreover $\ph'_n$ extends  to a cocycle.
Since $k>l$, the terms $\psd{1,0,xl}3$, which arise from our application of the reduction
formulas to $\ph'_x$ are not of the type that we use to extend $d_{k,l}$, because
$k>l$, and we should convert these terms to ones with lower index 2. Thus terms of the form $\psd{1,0,l+x-1}3$
are leading terms of coboundaries, and those of the form $\psd{1,0,k+(x-1)l}2$ are coboundaries when $x\ne n$.
Thus, any extension of $\d_{k,l}$ is of the form
$$
d_{k,l,m}=\psd{1,1,0}3+\psd{1,0,l}3+\psd{1,0,k}2+\psd{1,0,m}2a,
$$
where  $xm$ is not of the form $(x-1)l+k$ if $x\ne n$.
Using the modified recursion formulas (\ref{recur1},\ref{recur2},\ref{recur3}) from section \ref{Remark 2},
and using only the first term in the infinite series in \refeq{recur2}, since the others will clearly
lead to larger upper right indices,  we find after reduction that there are terms with upper right indices of the form
$i(k+1)+jl+z(m+1)+k-l$ and $i(k+1)+jl+z(m+1)+m-l$ for any triple $(i,j,z)$ of integers such that $2i+j+2z=n$
as long as $i\ge 1-n$, $0\le j\le n$, $z\ge 0$. The ones where $z=0$ correspond to the
terms from the unmodified reduction formulas.  These terms are handled by the same type of reasoning as in
the previous subcase 1.  The next lowest degree term is of type
$\psd{1,0,m+(n-1)l}2$, and its coefficient can be shown to equal
\begin{equation}
\frac{-b((n-1)l+1)((n+1)l-(m+1))}{nl}.
\end{equation}
This coefficient would vanish only if $m+1=(n+1)l$, which is impossible, since $m>k$. Thus our extended
codifferential picks up an additional leading coboundary term $\psd{1,0,m+(n-1)l}2$, and in extensions of
$\d_{k,l,m}$ we avoid adding terms of that type, as well as those we excluded before.

\centerline{\emph{ Subcase 2: $k+1$ is not a multiple of $l$}}

In this case,  all the even cocycles have been eliminated and $\psd{1,0,ml}3$ is the leading term of a coboundary
for all $m>0$. Thus, up to equivalence, and the action of the group $U_{k,l}$, extensions of $d_{k,l}$ are of the form
\begin{equation}
d_{k,l,e}=\psd{1,1,0}3+\psd{1,0,l}3+\psd{1,0,k}2+\sum_{v=k+1}^\infty\psd{1,0,v}2a_v,
\end{equation}
where $v\ne k+xl$ for any positive integer $x$.

\subsubsection{Extensions of $\dstar$ when $k+1=2l$}
Let
\begin{equation}
d_e=\psd{1,1,0}3+\psd{1,0,l}3.
\end{equation}
Recall from section \ref{Remark 2} that
$\ph_1(1-l)+\ph'_1l$ is the only non trivial $D_e$-cohomology class of order 1.
In general,  its $D_{k,l}$-coboundary is $-\psd{1,0,k}2a(k+1-2l)$, so that usually
$\psd{1,0,k}2$ represents a $D_{k,l}$.  But in the present case,  we see that $\ph_1(1-l)+\ph'_1l$ is
a $D_{k,l}$-cocycle, and $\psd{1,0,k}2$ is not a coboundary.
This is consistent with the observation that we cannot assume that the
constant $a$ appearing in $d_{k,l}$ is 1. Thus we begin with an extension of $d_e$ of the form
\begin{equation}
d_{k,l}=\psd{1,1,0}3+\psd{1,0,l}3+\psd{1,0,2l-1}2a.
\end{equation}

It is easy to see that every term which occurs in the reduced form of the coboundary of $\ph'_n$ has type
$\psd{1,0,nl}3$.
The coefficient occurring in the reduced form of $D_{k,l}(\ph'_n)$, for $n>1$ is
\begin{equation}
\frac{C(n)(a(n+1)^2+nl)\prod_{2p+1<n} (a(n-2p-1)^2-p(n-p-1)l)}{l^{[n/2]}},
\end{equation}
where $C(n)$ is a nonzero constant, depending only on $n$.
As a consequence, except when $n=1$, the coefficient of this term does not vanish for generic values of $a$.
In particular, when
$a=1$, this coefficient does not vanish for $n>1$.  On the other hand, for any value of $a$ which can be expressed
in the form $a=\tfrac{p(n-p-1)l}{(n-2p-1)^2}$, the coefficient vanishes, and moreover, if one of these coefficients
vanishes, then an infinite number of them do. In fact, if $p$ is the smallest value for which there is some $n$
such that $a=\tfrac{p(n-p-1)l}{(n-2p-1)^2}$, then $a=\tfrac{p'(n'-p'-1)l}{(n'-2p'-1)^2}$, iff $p|p'$ and
$(n-1)|(n'-1)$.
Therefore, for some very special values of $a$, we obtain a more
complicated phenomena than usual; especially, we see that there are an infinite number of values of $n$ for which
$\ph'_n$ extends to a nontrivial cocycle, so we have some nonzero even cohomology classes.

For generic values of $a$, we have only one even cohomology class, represented by $(1-l)\ph_1 +l\ph'_1$.
If we reduce the cochain $\psd{1,0,nl}3$ to lower index 2,
we get a nonzero multiple of $\psd{1,0,(n+1)l-1}2$, so
cocycles of the form $\psd{1,0,l+n-1}3$ and $\psd{1,0,(n+1)l-1}2$
are the leading terms of $D_{k,l}$-coboundaries. Thus, if we extend $d_{k,l}$, we would add a term of
the form $\psd{1,0,m}2b$ to it, where $m\ne (n+1)l-1$.

Whether or not $a$ is a very special value, if we consider any extension of the form
\begin{equation}
d_{k,l,m}=\psd{1,1,0}3+\psd{1,0,l}3+\psd{1,0,k}2a+\psd{1,0,m}2b,
\end{equation}
then by applying a diagonal automorphism, we can see that the coefficient $b$ can be taken to be 1.  In
fact,  $D_{k,l,m}(\ph_1(1-l)+\ph'_1l)=-\psd{1,0,m}2b(m+1-2l)$, so $\psd{1,0,m}2$ is a coboundary,
as we would expect, and there are no nontrivial 1-cocycles for the extended codifferential.

When $a$ is not a special value, then we restrict the terms which we add so that $m\ne(n+1)l-1$ for any
$n\ge 1$. When $a$ is special for the number $n$, consider what happens if we let $m=(n+1)l-1$. Then
for any value of $n'$ such that $(n-1)|(n'-1)$, we have  $D_{k,l}(\ph'_{n'})\sim 0$, and
\begin{multline}
D_{k,l,m}(\ph'_{n'})=\\D_{k,l}+\psd{1,n'-1,l(n+1)-1}2(n'+1-l(n+1))+\psd{1,n'-2,l(n+1)}3(n'-1),
\end{multline}
which is reducible to a multiple of $\psd{1,0,l(n+n')-1}2$. Every $n''$ such that $(n-1)|(n''-1)$ arises
exactly once in the set $n+n'-1$ as $n'$ ranges over the integers such that $(n-1)|(n'-1)$. Thus we see that
cochains of the form $\psd{1,0,(n''+1)l-1}2$,  where $(n-1)|(n''-1)$ also appear as
$D_{k,l,m}$-coboundaries of the former even $D_{k,l}$-cohomology class representatives.  Thus,  all the
even cohomology is killed, and we restrict any extension of $d_{k,l,m}$ to one in which the additional
terms are of type $\psd{1,0,x}2$, for $x$ not of the form $x=(n''+1)l-1$, for any $n''>0$.

We studied the reduction of $D_{k,l,m}(\ph'_{n'})$, and were not able to find a formula for the coefficient
of the leading term of type $\psd{1,0,(n+n')l-1}2$. On the other hand,  the terms in the reduced form
are of the types $\psd{1,0,p(n-1)+n'+1}2$, for $p=1\cdots n'$.  They are all of the form $\psd{1,0,(n''+1)l-1}2$
for some $n''$ such that $(n-1)|(n''-1)$, so that they don't appear as leading terms of $D_{k,l}$-coboundaries.
The highest degree term in $D_{k,l,m}(\ph'_{n'})$ is of type $\psd{1,0,(nn'+1)l-1}2$, and its coefficient is
not equal to zero.  Therefore, we know that $D_{k,l,m}(\ph'_{n'})$ gives rise to a new coboundary, whose
leading term is probably of type  $\psd{1,0,(n+n')l-1}2$.
If this remark is true, then we know exactly what are the types that must be excluded in the
extension process, that is, those which are of the form $\psd{1,0,x}2$, where $x=(y+1)l-1$ for any $y>n$.

Thus,  all the even cohomology is killed, and we restrict any extension
of $d_{k,l,m}$ to one in which the additional terms are of type
$\psd{1,0,x}2$, for $x$ not of the form $x=(n''+1)l-1$, for any
$n''>0$, if the leading coefficient is always  nonzero, and some other
pattern if the leading term does vanish.

\section{Conclusions} In this paper we have succeeded in classifying
all extensions of co\-differentials of degree less than or equal to 2
on a $2|1$-dimensional space.
Essentially, what we have been looking
at is the cohomology of \zt-graded Lie algebras, because a degree 2
codifferential determines precisely this structure on the space.

 Our purpose was not to give an
exhaustive classification of low dimensional \linf\ algebras, but to
explore the ideas that arise from the study of the classification and
extension problem in simple cases, to see what kind of interesting
phenomena can be observed.

As the reader can tell, even for simple examples of \linf\ algebras the
classification problem for extensions can be very intricate. The problem
of studying extensions seems to be far more complicated than constructing
a miniversal deformation.  Thus it is not surprising that the problem
of classifying nonequivalent extensions of an infinitesimal deformation to
a formal deformation should be difficult as well, since it should be at least
as complicated as the problem of classifying extensions of a \linf\ algebra
given by a codifferential of fixed degree. The problem of classification
of extensions of an infinitesimal deformation to a formal one is important,
and may well benefit from a study of examples similar to the ones we have
been addressing.



\providecommand{\bysame}{\leavevmode\hbox to3em{\hrulefill}\thinspace}
\providecommand{\MR}{\relax\ifhmode\unskip\space\fi MR }
\providecommand{\MRhref}[2]{%
  \href{http://www.ams.org/mathscinet-getitem?mr=#1}{#2}
}
\providecommand{\href}[2]{#2}

\end{document}